\documentclass[12pt,a4paper]{article}
\usepackage[top=16mm,left=12mm]{geometry}\textwidth18cm\textheight25.8cm
\usepackage{amsmath,amssymb}\usepackage{amsthm}
\usepackage{float}\floatstyle{plaintop}\restylefloat{table}
\usepackage[tableposition=top]{caption}
\usepackage{graphicx,pstricks,url,tikz}

\DeclareMathAlphabet{\mathdutch}{U}{dutchcal}{m}{n}
\SetMathAlphabet{\mathdutch}{bold}{U}{dutchcal}{b}{n}

\usepackage{drafts}
\usepackage{listings}\usepackage{paralist, longtable}
\theoremstyle{definition}
\newtheorem{theorem}{Theorem}[section]
\newtheorem{remark}[theorem]{Remark}
\newtheorem{lemma}[theorem]{Lemma}
\def\bexe{\begin{exercise}}\def\eexe{\eex\end{exercise}}
\def\bsol{\begin{solution}}\def\esol{\eex\end{solution}}
\def\bexa{\begin{example}}\def\eexa{\end{example}}
\def\brem{\begin{remark}}\def\erem{\end{remark}}
\def\bthm{\begin{theorem}}\def\ethm{\end{theorem}}
\def\blem{\begin{lemma}}\def\elem{\end{lemma}}
\def\bcor{\begin{corollary}}\def\ecor{\end{corollary}}
\def\bdefi{\begin{definition}}\def\edefi{\end{definition}}
\newcommand{\IDEA}{\textbf{Idea of the Proof.} }

\def\bmip{\begin{minipage}{\textwidth}}\def\emip{\end{minipage}}
\def\huga#1{\begin{gather} #1 \end{gather}}
\def\hugast#1{\begin{gather*} #1 \end{gather*}}
\def\hual#1{\begin{align} #1 \end{align}}
\def\hualst#1{\begin{align*} #1 \end{align*}}
\def\hueq#1{\begin{equation} #1 \end{equation}}
\def\hueqst#1{\begin{equation*} #1 \end{equation*}}
\newcommand{\btab}[2]{\begin{tabular}{#1}#2\end{tabular}}

\newcommand{\R}{{\mathbb R}}
\newcommand{\C}{{\mathbb C}}\newcommand{\N}{{\mathbb N}}
\newcommand{\Z}{{\mathbb Z}}

\def\CV{{\cal V}}\def\CA{{\cal A}}  
\def\CF{{\cal F}}\def\CG{{\cal G}}\def\CH{{\cal H}}\def\CL{{\cal L}}
\def\CO{{\cal O}}\def\CS{{\cal S}}
\def\CM{{\cal M}}
 \def\CW{{\cal W}}
\def\Cm{{\mathdutch {m}}}\def\Cr{{\mathdutch {r}}}

\def\setm{\setminus}

\def\ga{\gamma}\def\om{\omega}
\def\noi{\noindent}\def\ds{\displaystyle}
\def\vt{\vartheta}\def\pa{{\partial}}\def\lam{\lambda}
\newcommand{\bi}{\begin{itemize}}\newcommand{\ei}{\end{itemize}}
\newcommand{\ben}{\begin{enumerate}}\newcommand{\een}{\end{enumerate}}
\newcommand{\bce}{\begin{center}}\newcommand{\ece}{\end{center}}
\newcommand{\bci}{\begin{compactitem}}\newcommand{\eci}{\end{compactitem}}
\newcommand{\bcen}{\begin{compactenum}}\newcommand{\ecen}{\end{compactenum}}
\newcommand{\bcena}{\begin{compactenum}[(a)]}
\newcommand{\reff}[1]{(\ref{#1})}

\newcommand{\spr}[1]{\left\langle #1 \right\rangle}
\newcommand{\hs}[1]{{\hspace{#1}}}\newcommand{\vs}[1]{{\vspace{#1}}}

\def\eps{\varepsilon}\def\epsi{\epsilon}
\def\ra{\rightarrow}

\newcommand{\barr}{\begin{array}}\newcommand{\earr}{\end{array}}
\newcommand{\bpm}{\begin{pmatrix}}\newcommand{\epm}{\end{pmatrix}}
\newcommand{\bsm}{\left(\begin{smallmatrix}}
\newcommand{\esm}{\end{smallmatrix}\right)}
\newcommand{\ba}{\begin{array}}\newcommand{\ea}{\end{array}}
\def\dd{\, {\rm d}}\def\ri{{\rm i}}

\def\cc{{\rm c.c.}}\def\res{{\rm Res}}\def\er{{\rm e}}
\def\re{{\rm Re}}
\def\om{\omega}
\def\hot{{\rm h.o.t}}\def\ddt{\frac{\rm d}{{\rm d}t}}
\def\bphi{\varphi}\def\del{\delta}

\def\eex{\hfill\mbox{$\rfloor$}}

\def\Del{\Delta}

\def\sig{\sigma}
\def\al{\alpha}

\def\Ga{\Gamma}\def\Lam{\Lambda}
\def\kap{\kappa}

\def\bd{\begin{displaymath}} \def\ed{\end{displaymath}}
\def\ba{\begin{array}} \def\ea{\end{array}}
  
\def\eps{\varepsilon}

\def\pdep{{\tt pde2path}}

\newlength{\tew}

\def\ig{\includegraphics}

\renewcommand{\arraystretch}{1.05}\renewcommand{\baselinestretch}{1.0}
\def\medskip{}\def\bigskip{}
\usepackage{setspace}
\def\huc#1{{\magenta #1}} 
\def\amc#1{{\blue #1}} 
\def\huc#1{{#1}}\def\amc#1{{#1}} 
\def\stt{\small \tt}\def\sm{\small}\def\rb{\raisebox}

\usepackage{hyperref}
\def\taskip{\renewcommand{\arraystretch}{1}\renewcommand{\baselinestretch}{1}}
\def\teskip{\renewcommand{\arraystretch}{1.1}\renewcommand{\baselinestretch}{1.1}}

\def\hutab#1{\taskip\begin{\table}#1\end{table}\teskip}

\def\mS{\mathbb{S}}\def\MS{\mS}
\def\bsub{\begin{subequations}}\def\esub{\end{subequations}}
\def\qp{q_\text{phase}}\def\ddsi{\dd\sig}\def\Tau{\Upsilon}
\allowdisplaybreaks

\begin{document}
\text{}\vspace{10mm} 
\bce\Large
Breathing and moving vesicles 
in a geometric mechanochemical model\\[4mm]
\normalsize 
Alexander Meiners, Hannes Uecker\\[2mm]
\footnotesize
Institut f\"ur Mathematik, Universit\"at Oldenburg, D26111 Oldenburg, 
alexander.meiners@uni-oldenburg.de, hannes.uecker@uni-oldenburg.de\\[3mm]
\normalsize
\today
\ece

\begin{abstract} 
We consider a geometric mechanochemical model of vesicles which 
couples the Helfrich flow for 
the shape of a lipid bilayer vesicle membrane $X$ 
with a reaction-diffusion equation for a single ``morphogen'' $\phi$ on $X$.  
The Helfrich flow is the $L^2$ gradient flow of the 
 elastic bending energy $\CH(X)=\int_X (H-c_0)^2 dS$ of $X$, 
typically supplemented by area or volume constraints, or 
both.  The morphogen $\phi$ adsorbs/desorbs at places of high/low mean 
curvature $H$, i.e., the kinetics of $\phi$ depend on $H$, 
and conversely $\phi$ modifies the 
spontaneous curvature $c_0$ on $X$. 
The flow is no longer gradient, and hence allows for more complicated dynamics, 
including time 
periodic orbits, e.g., "breathing and moving" vesicle shapes. 
We show how to compute bifurcation diagrams for such 
solution branches via numerical continuation and bifurcation methods. We mostly 
focus on ``planar'' vesicles (1D closed curves)  
but also give an outlook on 3D vesicles (2D closed membranes). 
\end{abstract}
\tableofcontents

\section{Introduction} 
\begin{figure}[h]
\ig[width=1.05\tew]{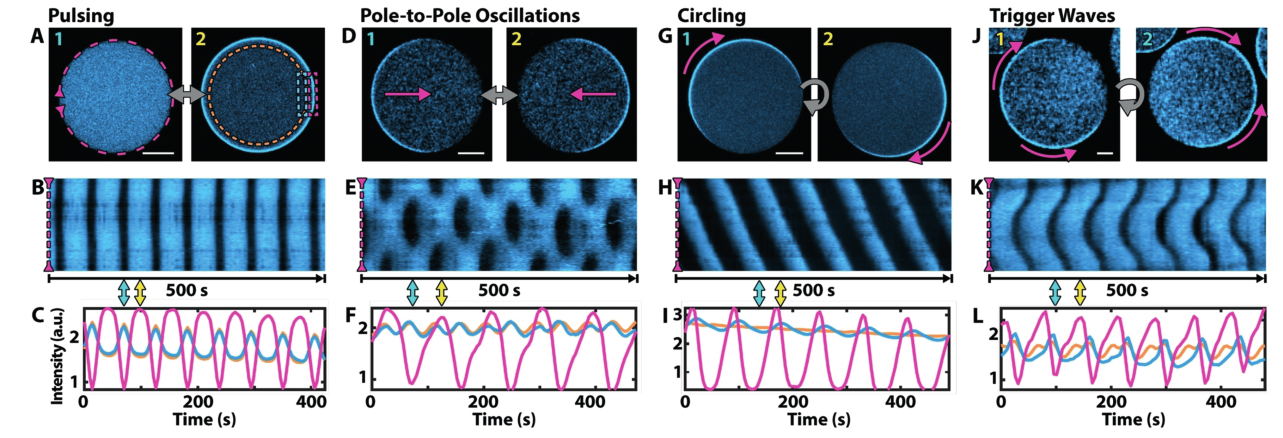}

\vs{-2mm}
\caption{{\sm Experiments from 
\cite{LRM18}. The top rows show snapshots of Min intensities, with the 
magenta arrows indicating the directions of waves. The 
kymographs and time series in the center and bottom rows are (normalized) 
intensities from the magenta and 
blue boxes in A, and the orange curve is the average intensity in the vesicle. }
\label{ef1}}
\end{figure}

Oscillations of biological membranes are an important 
building block of cell biology, and have been presented 
in the setting of giant unilamellar vesicles (GUVs) and other 
synthetic biomimetic vesicles in a number of recent studies 
\cite{TUSY14,LRM18,MSD18,SZHZ20,LDS20,LMAO21,TN22,L22,Nog25}. 
For instance, \cite{LRM18} reports experiments with 
GUVs that enclose a fluid with Min proteins that 
show some time periodic dynamics, and 
can change the spontaneous (or preferred) curvature $c_0$ of the membrane by 
adsorption to or desorption from the membrane. 
The membrane oscillations then found include ``pulsing'', where 
the MinD periodically in time but uniformly in space 
adsorbs/desorbs to/from the membrane,  
``Pole--to--Pole oscillations'' where MinD 
periodically switches between two spherical caps, ``circling'', and 
``trigger waves'', see Fig.\ref{ef1}. 
In all these, some area {\em and} volume changes of 
the vesicles seem involved, but no significant shape changes as the 
vesicles always stay essentially spherical. 
On the other hand, rather drastic shape (and topology) changes 
can also occur, for instance 
periodic closing and opening of necks of dumbbell shaped 
vesicles, and is explained 
by the phase diagrams of stable vesicle shapes from \cite{SBL90} 
depending on the spontaneous curvature $c_0$ of the Helfrich 
energy $\CH$ (see \reff{een2D}), namely that the periodic dumbbell 
dynamics follow the instantaneous minimizers of $\CH$ at the given $c_0(t)$. 

A classical geometric model of closed vesicle lipid bilayer membranes $X$ is 
the Canham--Helfrich energy \cite{H73} 
\hueq{ \label{een2D}
\CH (X, c_0)= \int_X (H-c_0)^2
  \,\dd S,
} 
where $H=\frac 1 2(\kap_1+\kap_2)$ 
is the mean curvature of $X$, i.e., the mean of the two principal 
curvatures on $X$, and the parameter $c_0\in\R$ 
is called the spontaneous curvature and depends on the 
properties of the lipid bilayers forming the membrane, and on the 
properties of the fluid in which the vesicles form. Stable vesicle shapes are 
obtained from minimizing $\CH$ under fixed enclosed volume $\CV$ and 
surface area  $\CA$, and the associated Euler--Lagrange equation 
is called the Helfrich equation. See, e.g., \cite{SBL90,SL95,sei97} for detailed 
2--parameter ($c_0$, and reduced volume $v=\CV/\CA$) 
phase diagrams of axisymmetric vesicles obtained from (numerical) solution 
of the Helfrich equation, and \cite[\S4.1]{geompap} 
for further continuation and bifurcation 
results, including non--axisymmetric shapes. 

The $L^2$ gradient flow associated to $\CH$ is called the Helfrich flow, and is 
studied from a mathematical perspective in, e.g., 
\cite{NY12,rupp24} for closed vesicles, and in \cite{mucyl} for 
periodic cylinders, which show wrinkling, pearling and coiling 
as primary instabilities of straight cylinders. 
However, as $\CH$ is non--increasing in these gradient flows, 
and in fact decreasing except in steady state, these models 
do not allow time-periodic solutions. 
To model time periodic behavior as in Fig.\ref{ef1}, a natural idea 
is to combine the gradient dynamics for $\CH$ with chemistry, 
yielding so called mechanochemical models, which 
come in at least two ways. 
\bci 
\item Bulk--surface coupling: In these, the dynamics of some chemicals in 
the fluid (the bulk) is coupled to the interaction 
of the chemicals with the membrane (the surface). 
The dynamics in the bulk may be of simple linear diffusion 
type. For instance, \cite{SG11} finds oscillations in a system of 
a fluid layer between two membranes which are graphs over a plane, where 
neither the bulk (linear diffusion) nor the membranes alone 
would generate oscillations. See also \cite{SGS07} for a somewhat 
simplified system without membrane shape changes, with a detailed analysis in \cite{GW16}, 
which is extended in \cite{PXP20}. 
On the other hand, in \cite{LRM18,CLSL21} it is argued that for the experiments 
shown in Fig.~\ref{ef1} the situation is rather the converse, 
at least in A--C: There is a 
time periodic reaction and fast diffusion of the Min proteins in the bulk, 
thus yielding an effective spatially homogeneous time-periodic forcing 
via $c_0(t)$ on the membrane. 
\item Reaction--diffusion (RD) systems on the membrane $X$: Here, 
some chemicals 
can adsorb/desorb from the fluid to $X$ and react and diffuse 
there, but the dynamics in the fluid are not modeled.  
The reaction on $X$, or in particular the adsorption/desorption, depends 
on the (mean) curvature $H$ of $X$, and vice versa the chemicals influence 
the spontaneous curvature $c_0$ on $X$. 
 A simple phenomenological model for this,  with just one abstract 
chemical species $\phi$ (a morphogen), is set up and analyzed 
numerically in \cite{MHMC13,MMRH13}, see also \cite{BMRM18}. 
Relatedly, in \cite{TN20,TN21}, the Helfrich flow for $X$ is coupled to
a Brusselator type RD system for two chemical species $v_1$ 
and $v_2$ on $X$ with dynamics depending on $H$, see also \cite{Nog25}. 
\eci 
Both classes of models come with different constraints, i.e., fixed 
membrane area, or fixed enclosed volume, or both. 
In any case, in both classes, the full dynamics are in general no longer 
of gradient type, and hence allow for time--periodic behavior.%
\footnote{There are also membrane models based on the Koiter 
elasticity model coupled to reaction--diffusion of chemicals in 
the surrounding (or enclosed) medium. For instance, in \cite{MSD18} 
a Koiter elastic energy is coupled to a FitzHugh--Nagumo model, while 
\cite{LDS20} 
considers disk--shaped Belousov--Zhabotinsky (BZ) gels in a 
BZ solution, and the periodic BZ reaction drives a swelling of the gel 
which alters its mechanical properties and hence yields periodic 
buckling, and a similar  
setup is considered in \cite{LMAO21}, with a Helfrich term included 
in the elastic energy. In all three models, time--periodic 
shapes changes of the membranes and periodic waves on the 
membranes are observed, which inter--alia may yield locomotion 
of the vesicles. 
}

A 1D version (without area or volume constraints) of 
coupling of Helfrich energy with branched $c_h$ vs bundled $c_v$ 
actin on membranes 
is treated in \cite{IMKYG13}, yielding a 4th order Helfrich type 
equation for the membrane height $h$ coupled to a two component 
RD system for $(c_h,c_v)$. See also \cite{G13}. 
Additionally, somewhat phenomenological 1D models without 
genuine geometry are treated in \cite{YFB22,BEGY23,HMY25}, i.e., 
reaction--diffusion models for F--actin, and 
active and inactive forms of GTPase, but without membrane mechanics. 
Under the crucial ingredient of conservation of the total GTPase mass, 
the models 
show bistability of ``wave--pinning'' (steady patterns), 
traveling waves, and standing oscillations, which 
are moreover related to experimental results for {\em D.~discoideum} cells. 

From the diverse possible mechanochemical models, and 
although the motivating experimental results from Fig.\ref{ef1} 
rather come from the bulk--surface class, here we analyze models 
based on \cite{MHMC13,MMRH13} via numerical 
bifurcation and continuation, 
following our previous work \cite{geompap,mucyl}. For closed 2D membranes, the 
model reads 
\bsub
\label{mm1}
\hual{
\pa_t X&=-[\Delta(H-c_0)+2H(H^2-K)+2c_0K-2c_0^2H-\lam H]N, \\
\pa_t \phi&=D\Delta\phi-\del\phi+g(H). 
}
In  the Helfrich flow (\ref{mm1}a) of $X$, again $H$ is the 
mean curvature (always wrt to the inner normal $N$ of $X$, i.e., $H=1/R$ 
for a sphere of radius $R$), $K$ is the Gaussian curvature, $\Delta$ is the 
Laplace--Beltrami operator on $X$, and the spontaneous curvature 
\huga{\label{c0f}
c_0=c_0(\phi)=\al+\beta\phi
} 
depends on $\phi$ via parameters $\al,\beta\in\R$.  The morphogen $\phi$ 
diffuses with rate $D>0$ on $X$ and decays with rate $\del>0$, 
and $H$ influences the dynamics of $\phi$ via the function 
\huga{
g(H)=\zeta \frac{f(H)}{\om+f(H)}, 
} 
with coupling constant $\zeta\in\R$ and parameter $\om>0$, and where 
$f(H)$ is a smoothed 
version of $\max(0,H-H_n)$ for some normalization $H_n>0$, e.g., 
$H_n=1/R$ for a sphere of radius $R$. The morphogen is 
thus ``generated'' (i.e., adsorbs) at places of higher $H$. 
Different from \cite{geompap,mucyl}, and as in \cite{MHMC13,MMRH13} 
we do not fix the enclosed volume which can adapt via osmosis, 
but the membrane $X$ is still inextensible, and 
$\lam$ is a Lagrange multiplier for the area constraint%
\footnote{Instead of the global area constraint $A(X)-A_0=0$, \cite{MHMC13,MMRH13} 
assume a local inextensibility and consequently introduce a field $\lam$ on $X$ and 
allow tangential motion of $X$, while here we restrict to normal motion only. 
The equivalence of both approaches for steady states of Helfrich flows of 
genus 0 surfaces (conformally equivalent to the sphere) is shown in \cite{DH15}. 
For the genuine dynamical 
problem \reff{mm1} the physical argument for using just a number $\lam$ is 
that the surface tension always relaxes so fast that it is effectively 
homogeneous and isotropic 
on $X$, cf.\cite{GNPS96}. Nevertheless, here this is a modeling choice/assumption.}
\huga{q(X)=A(X)-A_0=0.}
\esub

However, for numerical 
efficiency and simplicity, see Remark \ref{2drem}, we first consider 
1D membranes (curves) $\ga$. For these, 
the mean curvature $H$ of \reff{mm1} is replaced by the curvature $\kap$, 
hence $\CH (\ga, c_0)= \int_\ga (\kap-c_0)^2
  \,\dd s$, 
and the model reads 
\bsub\label{memo}
\hual{\label{helf1} 
\pa_t \ga &=-[\Delta (\kap-c_0)
+\frac 1 2 \kap(\kap+c_0)(\kap-c_0)-\lam\kap]\nu,\\
  \label{redi}\pa_t\bphi&=D\Delta \bphi 
-\del \bphi+\zeta\frac{f(\kap)}{\om+f(\kap)}, 
}
\esub 
where $\nu$ is the (inner) normal to the curve, 
together with the length constraint $q(\ga):=L(\ga)-L_0$, which 
determines the Lagrange multiplier $\lam$, where wlog $L_0=2\pi$. 

\begin{figure}[htp]
\btab{l}{{\sm (a)}\hs{26mm}{\sm (b)}\\
\ig[width=0.21\tew]{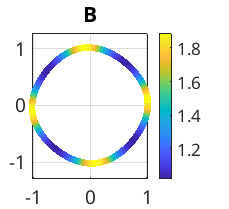}\qquad
\hs{-0mm}\ig[width=0.22\tew]{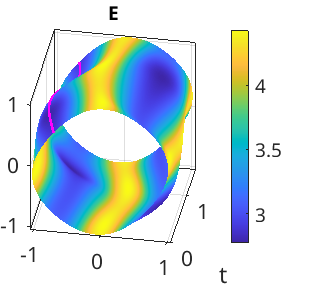}
\rb{12mm}{\btab{l}{\ig[width=0.2\tew]{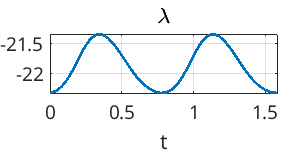}\\[-1mm]
\ig[width=0.2\tew]{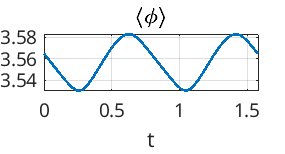}}}
\ig[width=0.28\tew]{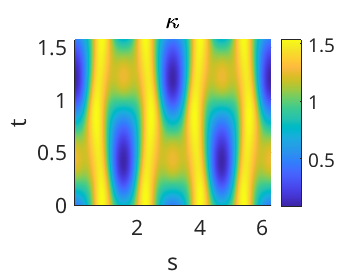}\\[-3mm]
{\sm (c)}\\[-2mm]
\hs{6mm}\ig[width=0.22\tew]{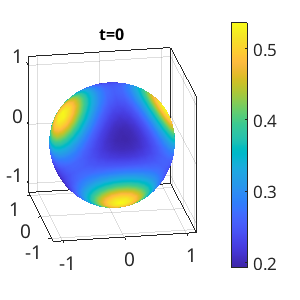}\ig[width=0.22\tew]{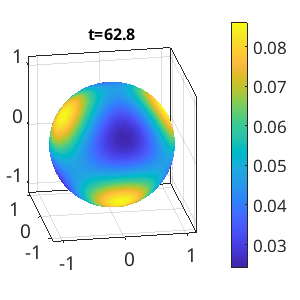}
\ig[width=0.28\tew]{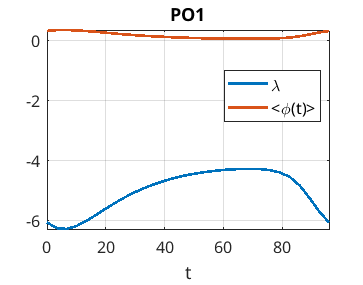}
}

\vs{-3mm}
\caption{{\sm (a) steady $D_4$ symmetric 1D vesicle, 
colored by $\phi$; $(\al,\beta,\del,D,\om){=}(3,1, 0.3,0.1,1)$, 
$\zeta\approx 5.67$. 
(b) breathing $D_4$--symmetric 1D vesicle (bifurcating from branch with B), 
$\zeta\approx 8.13$. 
Left: shape as a function of $t$; the red line indicates 
the arclength $s=0$ section. Middle: 
time series of $\lam$ and the average coverage $\spr{\phi(t)}:=
\frac 1 {2\pi}\int_0^{2\pi}\phi(s,t)\dd s$. Right: kymograph of $\kap$. 
(c) Sample of a breathing 2D  vesicle with  tetrahedral symmetry,  
$(\al,\beta,\del,D,\om)=(2,1,0.05,0.1,1)$, $\zeta\approx 3.45$. 
\label{f2a}}}
\end{figure}

Figure \ref{f2a}(a,b) gives a preview of 1D results, and (c) just 
one typical 2D time periodic orbit (PO). In both, 1D and 2D, we essentially fix 
the parameters $(\beta,\om,\del,D)$,   and for different choices of $\al$ 
compute bifurcation diagrams of steady states and POs with 
the coupling constant $\zeta$ as the primary active bifurcation parameter. 
In a nutshell, in 1D, the circle $\mS^1$ with $\phi\equiv 0$ 
is always a steady state 
for \reff{memo} (via $2\lam=c_0^2-\kap^2$) for the given parameters. 
This yields a ``trivial'' solution branch, but as 
the coupling constant $\zeta$ increases beyond some critical 
value, $\mS^1$ looses stability 
and first some nontrivial branches of $D_m$ symmetrical shapes 
bifurcate, see (a) for a $D_4$ symmetric sample. Then, POs 
may bifurcate from the $D_m$ symmetric shapes, sample in (b).  
The stability of these is difficult to assess from their Floquet 
multipliers, due to always present neutral translational and rotational 
modes, and hence we supplement the bifurcation analysis 
by some direct numerical simulation (DNS), see also Remark \ref{2drem}b). 
In these, the $D_4$ symmetric PO E (with kymographs like in Fig.\ref{ef1}E) 
turns out to be metastable in DNS in the labframe, and after some 
transient the solutions turns into ``trigger waves'' like in Fig.\ref{ef1}K, 
in with the vesicle additionally starts to tumble. 
However, if we integrate \reff{memo} in a comoving frame, i.e., 
trade spatial translations for positional Lagrange multipliers, see Remark \ref{2drem}b), then the PO E is stable. 
Similar behavior is presented for 2D vesicles in \S\ref{3dsec}, 
cf.~Fig.~\ref{f2a}(c). 

Importantly, while Fig.\ref{f2a} only 
previews samples (steady states and POs), in \S\ref{1Dsec} and 
\S\ref{3dsec} we give full bifurcation diagrams of branches of such states. 

\brem\label{2drem}{\rm 
a) The numerics are naturally much less expensive in 1D than in 2D; 
this is in particular important for POs, i.e., (1+1)D vs (2+1)D. 
Equally important, the 1D problem is much easier from 
a symmetry point of view. In 1D, the symmetry group of the 
primary bifurcations from the circle $\mS^1$ is $O(2)$, the linearization 
can be analyzed by Fourier modes, and the kernels 
are generically 2-dimensional. In 2D, we have $O(3)$ equivariant 
bifurcations from the sphere, with high--dimensional kernels given by 
spherical harmonics, cf.\cite{Matt04}. 

b) Due to the length constraint $q(\ga)=0$, 
\reff{memo} (resp.\,\reff{mm1} with area constraint $q(X)=0$ in 2D) 
always is a (partial) differential algebraic 
equation DAEs of (differentiation) index 2; moreover for continuation 
of steady states and POs, and also for DNS, we shall need additional 
positional constraints (aka phase conditions), yielding further 
index 2 components for the DAE system. In a simple setting, we then 
use an implicit Euler scheme for the DNS of these DAEs, which however 
we also cross--check with DNS based on the high--order and adaptive 
code {\tt RADAU} \cite{hlr89}%
\footnote{with code available at \url{https://www.unige.ch/~hairer/software.html}}, 
yielding essentially the same results. See also App.\ref{timevo}. 
\eex 
}
\erem 

\brem\label{mmrem1}{\rm  a) Given a base manifold $X_0$, and setting 
$X=X_0+u N$ with normal $N$ and normal displacement $u$, we can write 
(\ref{mm1}a), and the analogous version in 1D, as 4th order equations  
for $u$. In \cite{MMRH13,MHMC13} 
this mechanochemical model is then related to 2--component 
activator--inhibitor RD systems as follows: $u$ is a fast long--range 
(4th order) inhibitor, while $\phi$ (with typically small $D>0$) 
is a slow short-range activator. In this sense, one chemical species 
(the morphogen $\phi$) coupled to curvature is sufficient to obtain 
(tunable, see b)) 
pattern formation, and no genuine RD {\em systems} on $X$ are needed. 
However, this is also related to the important length (1D) resp.\,area (2D) 
constraints: Without these, i.e., without $\lam$ in \reff{helf1}, 
circular vesicles can just grow in radius and no patterns develop, 
and the same holds in 2D for growing spheres. The fast--slow point of view 
(small $D>0$, compared to fast 
relaxation of the membrane) is further pursued in the setting of \reff{memo} 
in \cite{NDV25} as fast--slow dynamics in space to 
study pattern formation  via geometric 
singular perturbation theory.

b) By ``tuneable'' in a) we mean 
that by adapting the chemical parameters $\beta,D,\delta$ and $\om$ 
(with $\zeta$ used as the continuation parameter) we can choose the 
wave-number $m{\in}\N$ (in 2D the spherical harmonic number $l$) at which 
under continuation in $\zeta$ the circle (the sphere) first becomes unstable, 
where we then expect bifurcation of stable patterns with wave number $m$ 
(spherical harmonic number $l$). In the 
purely mechanical Helfrich model, the first bifurcation is 
always with $m=2$ (resp.~$l=2$ in 2D, i.e., prolates and oblates). 
}
\eex\erem 

The rest of the paper is structured as follows. In \S\ref{1Dsec} 
we explain the 1D model \reff{memo}, 
including the precise form of $f(\kap)$, and discuss the linearization 
at the circle. We then present detailed BDs for this 1D model, 
and in \S\ref{3dsec} similar results in 2D, and in \S\ref{dsec} we end 
with a discussion and outlook. 
In App.\ref{lincal} we collect some results from differential 
geometry, discuss the setting of the phase conditions, 
and briefly the amplitude formalism 
to describe the steady bifurcations from the circle; this 
reflects the particular choice how we approximate the 
non--smooth function 
$f(\kap)=\max(0,\kap-\kap_n)$ of the original model 
by the smoothed versions $f_\chi(\kap)$ (see \reff{smax}). 
The numerical methods, in particular the extension of the {\tt Xcont} 
setting of \pdep\ to PO problems, 
the general mesh--handling strategies, and the DNS methods for the DAEs are 
discussed in \S\ref{algo}.  The main bottleneck is the computation 
of POs, and in particular in 2D we use rather coarse meshes. 
The numerical setting (mesh sizes, tolerances) is summarized in 
\S\ref{ncsec}, and we give some examples of adaptive mesh refinement 
in space and time, illustrating that the qualitative features are well 
resolved, although in 2D some quantitative changes occur under mesh refinement.

\begin{table}[h]
{\small \caption{Acronyms.}\label{tab1}
\bce
\vs{-2mm}
\begin{tabular}{|l|l|l|l|}   \hline
\bf Acronym & \bf Meaning&\bf Acronym & \bf Meaning \\ 
\hline 
AE&amplitude equation&BD&bifurcation diagram\\
BP, FP, HP&branch/fold/Hopf point&DAE&differential algebraic equation\\
DNS&direct numerical simulation&cDNS&constrained DNS\\
IC&initial condition&PO&(time-)periodic orbit\\
COM&center of mass&RD&reaction diffusion\\
\hline
\end{tabular}
\ece
}
\end{table}

\noi{\bf Notation.} Table \ref{tab1} summarizes the 
acronyms used in the paper. Additionally, by $D_m$ we denote
the dihedral group of order $2m$, i.e.,  
the symmetry group of a regular $m$--gon, generated by rotation with angle 
$2\pi/m$ and one reflection, 
important for the symmetry considerations of bifurcations 
in the 1D model. \\[2mm]
\noi{\bf Acknowledgement.} We thank the DFG for support under 
UE 60/6-1.

\section{The 1D model}\label{1Dsec}
For the 1D version \reff{memo} of \reff{mm1}, let 
\hueq{ \label{een}\CH (\ga, c_0)= \int_\ga (\kap-c_0)^2
  \,\dd s,
} 
be the bending energy of a (closed) curve $\ga\in\R^2$ with curvature $\kap$, 
where again $c_0\in\R$ is a spontaneous curvature, which we take 
to depend on a morphogen $\phi$ via \reff{c0f}, i.e., 
$c_0{=}c_0(\phi){:=}\al{+}\beta\phi$ 
with parameters $\al{\in}\R$ and $\beta>0$. 
The gradient flow of \reff{een} with respect to the $L^2$ inner product,
$\spr{f,g}_{L^2}:=\int_\ga \spr{f,g}_{\R^2}\dd s$, under the length 
constraint $L(\ga)=L_0$, and 
coupled with the morphogen dynamics yields \reff{memo}, 
which below we express in local coordinates. 

The original model uses $f(\kap)=\max(0,\kap-\kap_0)$ with 
$\kap_0$ a reference curvature, 
e.g., $\kap_0=1/R$ for a circle of radius $R$. However, this is 
not differentiable at
$\kap=\kap_0$, and therefore we replace it with the smooth approximation 
\amc{$\ds f(\kap):=\frac 1 \chi \log\left(e^{\chi(\kap-\kap_0)}+1\right)$} 
which converges to $\max\{\kap-\kap_0,0\}$ as $\chi\ra \infty$, uniformly 
for $\kappa\in\R$. For convenience 
(i.e., to have the trivial branch $(\kap,0)=(\kap_0,0)$ independent of 
further parameters) 
we  subtract $\log(2)/\chi$ from $f(\kap)$, i.e., finally choose 
\hueq{\label{smax}\amc{f_\chi (\kap):=\frac 1 \chi\left(
\log\left(e^{\chi(\kap-\kap_0)}+1\right)-\log(2)\right),}}
and generally fix $\chi=50$, see Fig.\ref{ampf1} in \S\ref{aesec} for sketches. 
We then have eight parameters 
$(L,\kap_0,\al,\beta,D,\\ \del,\zeta,\om)$, but by scaling 
we can reduce to $L=2\pi$, $\kap_0=1$, i.e.: 
The elastic energy satisfies
\amc{$\CH(\eps_\ga \ga,c_0)=\eps_\ga^{-1} \CH(\ga,\eps_\ga^{-1}c_0)$} for any
$\eps_\ga>0$; hence \reff{helf1} shares this property. In particular,
the case $c_0=0$ (Willmore energy) is scale-invariant. 
Thus, we can fix the scale by choosing $L=2\pi, \kap_0=1$, and 
adapting $\al,\beta$. Here 
we restrict ourselves to  the case of fixed $(\beta,\om)=(1,1)$, 
and $(D,\del)=(0.1,0.3)$, as the dynamics of $\phi$ is considered to 
be slower than that of $\ga$, and will look at different 
values of $\al$, while the coupling constant $\zeta$ 
will be our primary continuation parameter.

\brem\label{smrem}{\rm The smoothness of $f{=}f_\chi$ from \reff{smax} is needed 
for standard local stability and bifurcation analysis of $\mS^1$, and similarly 
at the sphere in 2D.  
We have $f(\kap_0)=0$ and $f'(\kap_0){=}1/2$ independent of $\chi$, 
but $f''(\kap_0){=}\chi/4$ and $f^{(n)}(\kap_0){=}\CO(\chi^{n-1})$.  
Thus, the linear stability of $\mS^1$ (with $\kap_0{=}1$) does 
not depend on $\chi$, but the higher order terms in 
the local bifurcation analysis for \reff{memo} at 
BPs from $\mS^1$ strongly depend on $\chi$. In particular, 
the validity range of the amplitude equations for the bifurcations 
shrinks with increasing $\chi$,  see \S\ref{lincal}. Thus, 
the choice of $\chi$ is a modeling decision, and 
we took $\chi{=}50$ to be close to the original model, and $f_\chi(\kap){<}0$ 
for $\kap{<}\kap_0$ can be seen as an adjustment of the baseline morphogen 
level at the circle, which in turn can be compensated by adjusting 
$\al,\beta$ in \reff{c0f}.
}
\eex 
\erem

\subsection{Local coordinates, local existence, 
and linear stability analysis at the circle} 
\label{locstab}
Let $\ga :[0,2\pi)\ra \R^2$ be a smooth, parameterized, inward-oriented 
closed curve, and let $\tau$ and $\nu$ denote 
the unit tangent and normal vector fields along $\ga$, defined by
\hueq{ \tau=\frac{\ga_s}{\|\ga_s\|_2}, \quad \nu= \frac{\tau_s}{\|\tau_s\|_2},}
where $\ga_s=\pa_s \ga$, where by smoothness we can assume $\ga_s\ne 0$ throughout. 
We denote the metric by $g=\|\ga_s\|_2^2$, and the 
Laplace--Beltrami operator on $\ga$ by $\Delta$, with 
\amc{$\Delta u=g^{-1/2}\pa_s\left(g^{-1/2}\pa_s u\right)$}, and then have 
\hueq{\kap=\spr{\Delta \ga, \nu}}
for the curvature $\kap$. 

To evaluate \reff{helf1}, \cite{MMRH13} 
considers variations of the form $\ga+u\,N+\psi\, T$, 
involving both the normal component $u$ and the tangential component $\psi$. 
Here, wlog (see \cite{KPP17}) we restrict to normal variations. 
If we write $\ga_\eta=\ga+\eta u\nu$, and with a slight abuse of 
notation let $L(\eta)=L(\ga_\eta)$, then $L'(0)=-\int \kap u \dd s$. 
Therefore, \reff{helf1} becomes 
\hueq{\label{nflow} \spr{\pa_t \ga,\nu}=-\bigg(\pa_u\CH -\lam(u) \kap\bigg),
\quad \pa_u\CH=\Delta (\kap-c_0)
+\frac 1 2 \kap(\kap+c_0)(\kap-c_0), 
}
and where $\lam(u){=-}(\int_\ga \kap^2\dd s)^{-1}\int_\ga \kap\pa_\ga \CH \dd s$ 
 is a non-local term  obtained from $0=\ddt q(\ga)=\pa_\ga L(\ga)[\pa_t\ga]$.  
The notation $\lam(u)$ in \reff{nflow} thus means that 
we consider the Helfrich flow projected on the constraint $L(\ga)=L_0$. 

The trivial steady state solution is the circle $\mS^1$, with 
 $2 \lam=c_0^2-\kap^2$ for any $c_0\in \R$. 
We fix the 
arc-length parametrization of $\mS^1$ so that $\kap=1$, i.e., 
$s\mapsto (\cos(s),\sin(s))^T$ and $\nu_0=-(\cos(s),\sin(s))^T$, and 
let 
\hueq{\ga=\mS^1- u \nu_0,
}
\textit{where we chose $-u$ instead of $+u$} as the formulas become 
slightly nicer. 
We rewrite \reff{nflow} for a normal variation $u$ of $\mS^1$ with the metric
\hueq{ g(u)=(1+u)^2+u_s^2,
}
and the curvature
\amc{\hueq{\kap(u)=\frac 1 {g^{1/2}(u)}+\frac {2u_s^2-(1+u)} {g^{3/2}(u)}u_{ss}, 
}}
and obtain 
\hueq{\label{locmemo} 
\bpm \frac {1+u} {g^{1/2}(u)} & 0\\ 0& 1\epm 
\bpm u_t\\ \phi_t\epm= \CG(u,\phi)=
\bpm - \left(\Delta (\kap-c_0)+\frac 1 2 \kap\big(\kap+c_0\big)\big(\kap-c_0\big) -\lam(u) \kap\right) \\
D\Delta \phi -\del \phi+\zeta f_\chi(\kap)(\om+f_\chi(\kap))^{-1}
\epm.
}
Following \cite{mucyl}, we define \def\epsi{\epsilon}
\hueqst{V_\epsi^\theta=H^\theta(\mS^1)\cap \{u:\, L(\mS^1-u\nu)=2\pi\} 
\cap \{u\ge -1+\epsi\},
} 
for $\theta>5/2$ (not a half-integer) and $0<\epsi<1$. 
Then $V^\theta_\epsi$ is non-empty (spanned by Fourier modes with 
wave number $m\in\Z\setm\{0\}$ with $0$ excluded by the length 
constraint), 
and \reff{locmemo} is parabolic and possesses the maximal regularity 
property on $V^\theta_\eps\otimes H^1(\mS^1)$ as a closed subset of 
$L^2(\mS^1)\otimes L^2(\mS^1)$. From this we obtain local 
existence for \reff{locmemo} for initial conditions 
$(u_0,\phi_0)\in(V_\eps\otimes H^1)$, cf.\,\cite{mucyl}, 
with simple modifications to account for the $\pa_t\phi$ equation.%
\footnote{This works fully analogous near the sphere in 2D, with the 
only difference that $\theta>3$ is needed then.}
However, we have no general results on global existence. 

For the numerical algorithms, it is helpful to also implement the Jacobian 
$\pa_{(u,\phi)}\CG$ at a general $\ga$, see \S\ref{lincal} 
for the derivation, 
but here we first explicitly linearize at $\mS^1$, i.e., 
at the homogeneous steady state $\bpm u_h,\phi_h\epm= \bpm 0,0\epm$, 
with $\lam(0)= \frac 1 2(c_0^2-1)$.
For normal variations of $\mS^1$, the linearization 
of \reff{memo}  reads 
\hueq{\label{linflow} \bpm 1 & 0\\ 0 &1\epm U_t = 
\CL(\pa_s,\Lam) U, 
}
where $\Lam$ stands for the parameters $(D,\del,\beta,\om,\al,\zeta)$. 
This is a linear parabolic system with constant coefficients and hence 
has solutions of the form 
\hueq{ 
U(s,t)=\exp(t\mu_m+\ri ms) \Phi_m\,
}
for any wave numbers $m\in \Z\setm\{0\}$ (again $0$ excluded by the length 
constraint), where $(\mu_m,\Phi_m)\in \C\times\C^2$ is an eigenpair 
of the Fourier-transformed linearization 
\hueq{\label{Lf}
\CL(\ri m,\Lam)=\bpm
-(m-1)^2(m+1)^2&
\beta(m^2-\al)\\ 
\frac {\zeta }{2\om}(m^2 -1)&
-D m^2-\del
\epm.}
Due to $\int_\ga e^{ims}\dd s=0$ for $m\neq 0$ we have 
$\pa_u\lam(u)|_{u=0}=0$ such that the constraint is absent in \reff{Lf}. 
For $m=1$, we have two zero eigenvalues related to rigid body motions, 
which we will remove via phase conditions, see \S\ref{pcs}. 
Also note that $\pa_u\CG_1$ is independent of any parameter; 
hence, $\mS^1$ is stable without morphogen coupling. 
This is different from the work in \cite{mucyl}, where the presence 
of pressure destabilizes the circle (treated as a cylinder in 
\cite{mucyl}). See also \S\ref{zezero} for a review of this case. 

\begin{figure}[h]
\btab{lll}{
{\small (a)} & {\small (b)} &{\small (c)}\\
\raisebox{0mm}{
\hs{-3mm}\ig[width=0.39\tew]{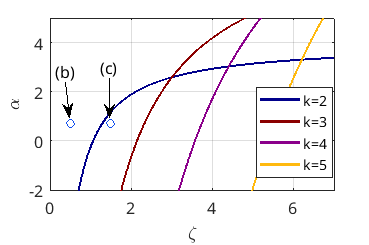}} &
\hs{-3mm}\ig[width=0.32\tew]{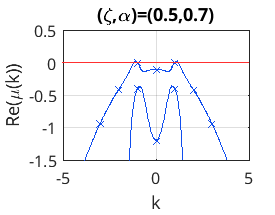}&
\hs{-3mm}\ig[width=0.32\tew]{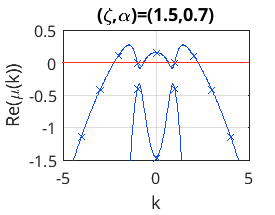}
}

\vs{-3mm}
\caption{\small{(a) Bifurcation curves obtained 
from \reff{zetc} for $m=2,3,4,5$ in the $\zeta,\al$---plane for 
fixed $(\beta,\om,D,\del)=(1,1,0.1,0.3)$. (b,c) Eigenvalues at points marked 
in (a). 
}  
\label{specplot}}
\end{figure}

By analyzing the dispersion relation $m\mapsto \mu_m(\Lam)$ 
for \reff{Lf}, we find that an easy way to control the wave number of 
the primary instability is by varying 
$\al$. Namely, given a wave number $m$, solving 
$\mu(m,\Lam)=0$ for the critical value of our bifurcation parameter $\zeta$  
we obtain 
 \hueq{\label{zetc}
\zeta_c(m) =\frac{2\om(m - 1)(m + 1)(Dm^2 + \del)}{(m^2 - \al)\beta}, 
}
which for instance gives the plots in Fig.\ref{specplot}(a), with 
subsequent plots of the eigenvalues at the marked points. 
Of course, there are other possibilities to manipulate the 
$\zeta_c$ (in particular, $\zeta_c(m)$ is a \amc{increasing} function 
of $\om, D$, and $\del$), but \reff{zetc} and Fig.\ref{specplot}(a) show that, 
e.g., for 
fixed $(\al,\del)=(0.7,0.3)$  (and always 
$(\beta,\om,D)=(1,1,0.1)$ fixed), increasing $\zeta$ from $\zeta<\zeta_c$ 
the first bifurcation is to $m=2$, the second to $m=3$, and so on. 
On the other hand, for $\al=3$, the order is $m=3,2,4,\ldots$, 
and for instance for $\al>4$ the $m=2$ mode never becomes unstable.%
\footnote{An analogous analysis is possible at the 2D sphere, and shows 
how to ``tune'' the bifurcation order from the sphere. However, this 
becomes considerably more complicated, as the Fourier modes $e^{\ri ms}$ 
must be replaced by spherical harmonics, cf.\,Rem.\,\ref{2drem}(a).} 
Moreover, further analysis shows that no Hopf bifurcations 
can occur on $\mS^1$ and hence to search for POs we must first follow the bifurcating 
steady state branches and ``hope'' for Hopf bifurcations on these. 
Essentially, these are also the reasons why we choose positive $\al$: 
 For small 
$\al$ ($\al<0.5$, say) we find fewer HPs on the nontrivial branches, 
in particular no HPs on the $m=2$ branch for $\al<0$, and 
the rather large $\al=3$ was chosen to put the $m=3$ bifurcation 
before the $m=2$ bifurcation. 

\subsection{Numerical continuation}\label{num}
We now use numerical continuation 
to compute the steady state branches $\CW_m$ with wave number 
$m\ge 2$ bifurcating at $\zeta_c(m,\Lam)$ from \reff{zetc}, 
and secondary bifurcations from these, including branches 
of POs arising from Hopf bifurcations. 
Throughout we fix 
\huga{\label{bps}
\text{$(\beta,D,\om,\del)=(1,0.1,1,0.3)$}} 
and use $\zeta$ as the bifurcation parameter, for two cases 
$\al=0.7$ and $\al=3$. 
Additionally, for the sake of completeness, we first discuss the 
case  $\beta=0$. 
Here, the membrane decouples from $\phi$ and its dynamics are 
governed by the gradient Helfrich flow of the energy
\reff{een}, with fixed $c_0=\al$.  
For any given IC $\ga|_{t=0}$ with global flow, the $\om$--limit set 
then consists of a steady state shape $\ga_\infty$ with 
associated Lagrange multiplier $\lam_\infty$.  
The BD thus consists only of steady states, 
and as  $\phi$ follows $\kappa$ we for 
simplicity drop the $\phi$ equation completely. 

Subsequently, going back to the case with coupling, i.e., $\beta=1$, 
we explore the more global behavior of the $\CW_m$ branches, 
including their loss or gain of stability away from onset, 
and some secondary bifurcations of steady states, 
and of {\em relative} steady states, meaning rigidly 
drifting fixed shapes and $\phi$ coverages; these are computed 
as steady states with non--zero positional Lagrange multipliers 
$\tau_x$ and $\tau_y$, see \S\ref{pcs}.%
\footnote{In principle, there could also be rigidly rotating shapes, but such were 
not found, and instead we only find rotating relative POs, see below.} 
As our main objective 
we then look at POs bifurcating at HPs from the (relative) steady state branches,  
and observe two major classes of interesting POs:
breathers (breathing cells without net motion), and breathing {\em and} 
moving cells. 
In both, the length constraint restricts the amplitude of the POs, 
and the Lagrange multiplier $\lam$ for the length constraint must be treated as a dynamical 
variable, and similar we need Lagrange multipliers $\tau_x,\tau_y$ and $\rho$ for 
translation  and rotation, respectively, thus 
yielding comoving frames with periodic speeds. 
We call the POs 
in these comoving frames {\em relative POs}, and plot 
these in the comoving frame, augmented by $t\mapsto (\tau_x,\tau_y,\rho)(t)$ 
if these are non--zero, and by the paths of the centers of mass of the cells 
in the labframe (e.g., Fig.\ref{al07h1}(c) and Fig.\ref{al3c2}(b--d)). 
Additionally, we use DNS to study the stability of POs. 
From this, besides 
the different ordering of $m$ for the primary bifurcations, another 
difference between $\al=0.7$ and $\al=3$ is that for the former the 
POs found by Hopf bifurcation are all unstable, while for $\al=3$ 
we find (meta)stable breathing $\CW_4$ branches, bifurcating supercritically 
at a loss of stability of the steady $\CW_4$ branch.  

All numerics are based on \pdep\ \cite{p2pbook,p2phome}, in particular 
extending the {\stt Xcont} setting \cite{geompap,mucyl} to PO 
computations, see \S\ref{algo}. Software sources for and 
further documentation of the computations can be found at \cite{mcSI}. 

\subsubsection{Intermezzo, $\beta{=}0$: Destabilization via external pressure}
\label{zezero}
The analysis in \S\ref{locstab} shows that with $\phi$ ``switched off'',  
$\mS^1$ is stable for all $c_0\equiv \al$. One way of destabilizing $\mS^1$ in 
the purely mechanical model is by 
introducing an osmotic pressure $P$, which can also 
be seen as a Lagrange multiplier for an enclosed volume constraint 
(in 1D: enclosed area constraint). Hence, for $\beta=0$ we consider  
\hueq{\label{helfpre}\spr{\ga_t,\nu}=\CG_1(u)+P.}
Without loss of generality we can set $c_0=0$, as $c_0$ only appears as $c_0^2$ 
in front of $\kap$, and serves the same purpose as $\lam$; then 
$\mS^1$ is a steady state for any $P$, with $\lam=\frac 1 2 -P$. 

\begin{figure}[h]
\btab{llll}{
{\small (a)}\hs{45mm}{\small (b)}\\
\btab{l}{ 
\hs{-3mm}\ig[width=0.40\tew]{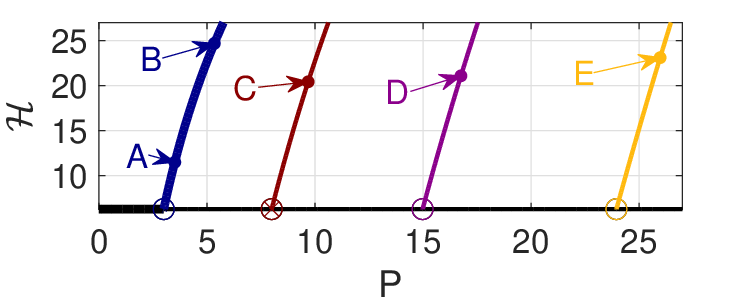}\\
\hs{-3mm}\ig[width=0.40\tew]{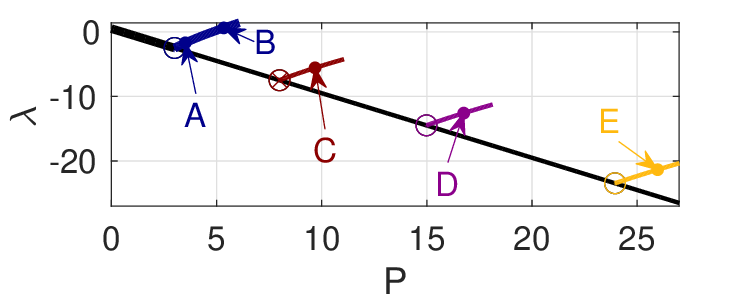}}
\raisebox{3mm}{
\btab{ll}{
\raisebox{2mm}{ \hs{-2mm}\ig[width=0.25\tew]{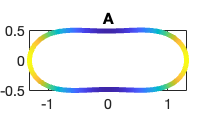}}
\raisebox{0mm}{ \hs{0mm}\ig[width=0.3\tew]{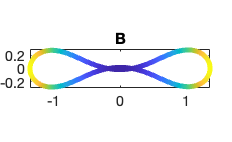}}\\
\raisebox{0mm}{ \hs{-2mm}\ig[width=0.18\tew]{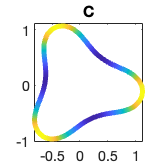}}
\raisebox{0mm}{ \hs{0mm}\ig[width=0.18\tew]{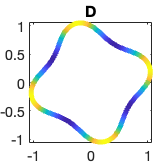}}
\raisebox{0mm}{ \hs{0mm}\ig[width=0.18\tew]{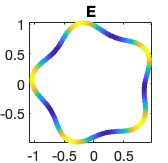}}
}
}

}
\vs{-3mm}
\caption{\small{(a) Basic bifurcation diagram for \reff{helfpre}, 
Helfrich energy $\CH$ (top) and Lagrange multiplier $\lam$ (bottom) over the pressure $P$. 
Thick lines indicate linearly stable branches. (b) sample solutions. }
\label{beta0fig}}
\end{figure}

Figure \ref{beta0fig} shows basic BDs  
(BPs as $\circ$, stability/instability via thick/thin lines)
and sample solutions. The black branch corresponds 
to $\mS^1$ with $\lam=\frac 1 2-P$ and constant area $\pi$. 
At $P=3$ the circle destabilizes in a 
pitchfork bifurcation to the $\CW_2$ wrinkling branch with 
decreasing area as $P$ increases. The branch 
remains stable even beyond the point of self intersection, reached 
shortly before sample B. 
The further primary bifurcating branches $\CW_3,\CW_4,\ldots$ 
stay unstable throughout, and like $\CW_2$ show self intersections 
at larger amplitudes. This is somewhat nonphysical, but 
clearly supported by the model.  
There are no secondary bifurcations on 
any $\CW_k$ branches up to self intersection, quite in contrast to the 2D case 
\cite{SBL90,geompap}, and also very different from the 1D case 
with morphogen, see below. Thus, by continuation and bifurcation we 
can only find $\mS^1$ and $\CW_2$ as stable branches, and this 
raises the question, if for any given $P$ 
we can find other stable solutions via 
(length $L=2\pi$ preserving) 
DNS. For this we perturbed a variety of (unstable) steady states 
from the $\CW_3,\CW_4,\ldots$ branches and used these as ICs 
for the Helfrich flow, and in all cases found convergence to 
the (possibly self intersecting) $\CW_2$ solution at the given length. 
This shows that the 1D Helfrich model is rather restricted 
in its (length preserving) dynamics. We are now ready 
to look at the mechanochemical model.

\subsubsection{$\al=0.7$} \label{alon}
\begin{figure}[h]
\centering
\btab{l}{
\btab{lll}{{\sm (a)}&{\sm (b)}&{\sm (c)} \\[-1mm]
\hs{0mm}\ig[width=0.43\tew]{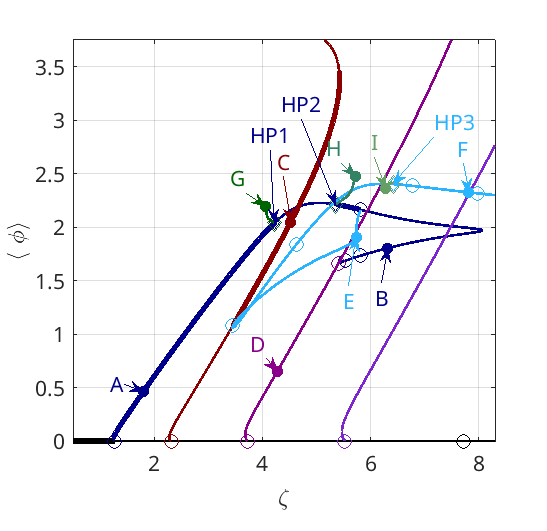}
&\hs{-5mm}\rb{3mm}{\ig[width=0.27\tew]{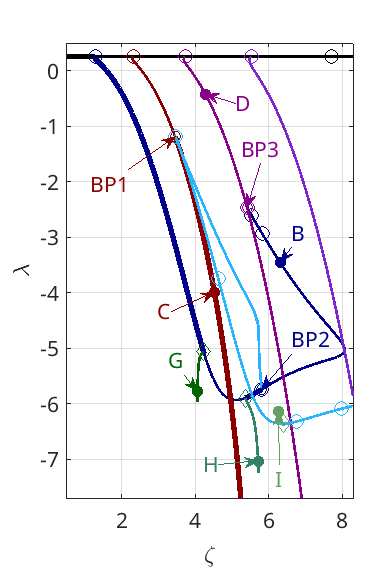}}
&\hs{-5mm}\rb{3mm}{\ig[width=0.27\tew]{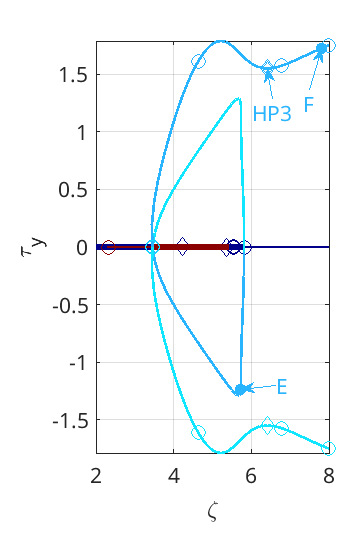}}
}\\[-6mm]
{\sm (d)}\\[-1mm]
\hs{0mm}\ig[width=0.18\tew]{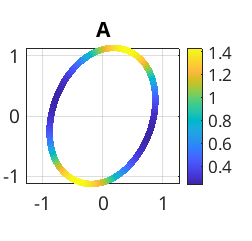}\hs{0mm}\ig[width=0.18\tew]{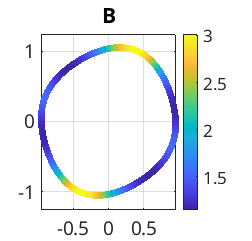}
\hs{-3mm}\ig[width=0.18\tew]{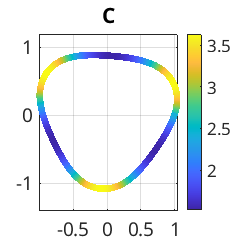}\hs{0mm}\ig[width=0.18\tew]{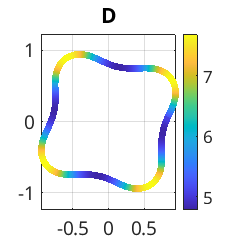}
\hs{-3mm}\ig[width=0.18\tew]{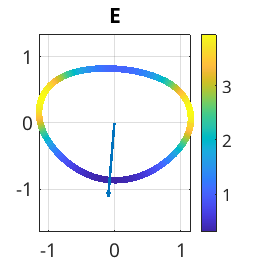}\hs{-3mm}\ig[width=0.18\tew]{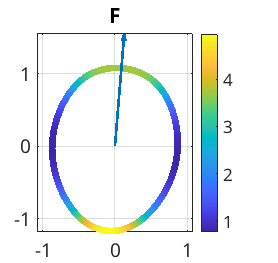}}\\[2mm]
{\sm 
\begin{tabular}{l|cccccccccc}  \hline 
&A&B&C&D&E&F&G&H&I\\
$\zeta$&1.8043&6.3136&4.5225&4.8466&5.7326&7.8153&4.0610&5.7137&6.2716\\
\hline
\end{tabular}}

\vs{0mm}
\caption{\small{(a) The basic bifurcation diagram (BD) for \reff{memo} 
with parameters from \reff{bps} and $\al=0.7$; $\spr{\phi}$ 
over $\zeta$. (b) same branches with $\lam$ over $\zeta$, where $\lam=\spr{\lam}$ 
for the (green) PO branches. (c) behavior of $\tau_y$ in and after the drift 
bifurcation(s) of the light blue branch. The very light blue branch in (c) corresponds 
to the ``other directions'', i.e., with the arrows in samples E and F flipped. 
(d) steady state samples; the shapes and arrows in E and F indicate the drift; 
see Fig.\ref{al07h1} for PO samples. }  
\label{al07f1}}
\end{figure}

\begin{figure}[h]
\btab{ll}{
\hs{-5mm}\btab{l|l|ll}{{\sm (a)}&{\sm (b)}&{\sm (c)}\\[-0mm]
\hs{0mm}\btab{l}{\ig[width=0.24\tew]{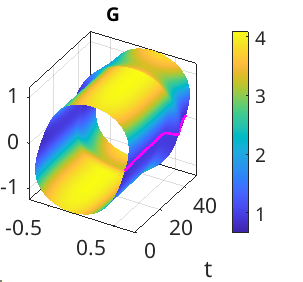}\\
\hs{0mm}\ig[width=0.22\tew]{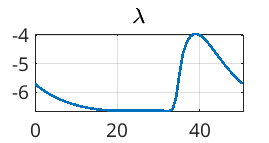}\\
\hs{0mm}\ig[width=0.22\tew]{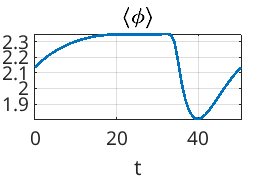}\\
\hs{-5mm}\ig[width=0.26\tew]{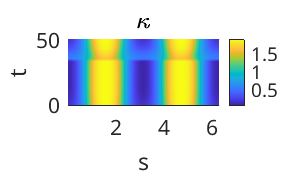}}
&\hs{-0mm}\btab{l}{\ig[width=0.24\tew]{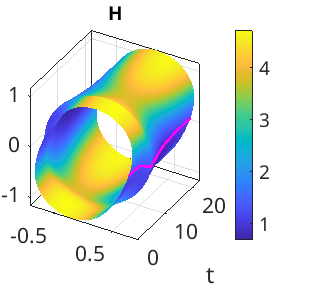}\\
\hs{0mm}\ig[width=0.22\tew]{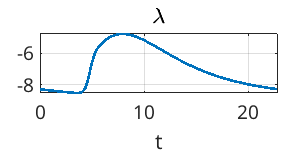}\\
\hs{0mm}\ig[width=0.22\tew]{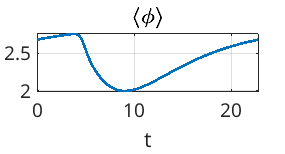}\\
\hs{0mm}\ig[width=0.26\tew]{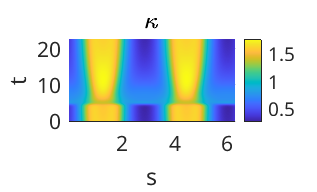}}
&\hs{-0mm}\btab{l}{\ig[width=0.24\tew]{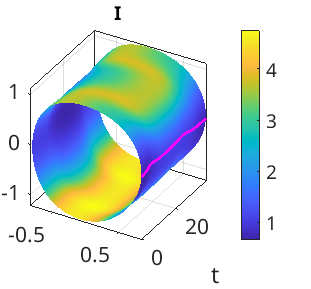}\\
\hs{0mm}\ig[width=0.22\tew]{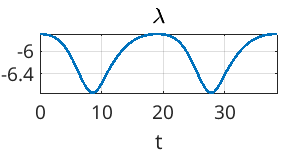}\\
\hs{0mm}\ig[width=0.22\tew]{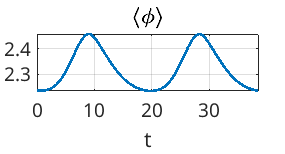}\\
\hs{-0mm}\ig[width=0.26\tew]{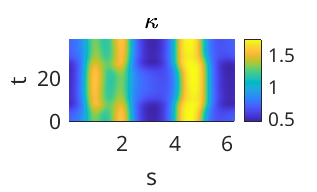}}
}
\hs{-12mm}\rb{10mm}{\btab{l}{\ig[width=0.21\tew]{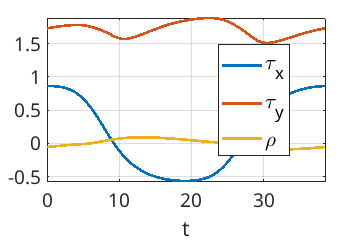}\\
\hs{-1mm}\ig[width=0.21\tew]{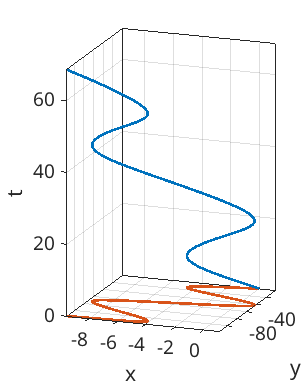}}}
}
\vs{-2mm}
\caption{{\sm PO samples G (a), H (b) and I (c) from the BD in Fig.\ref{al07f1}. 
Top to bottom: shape of $\ga$ along 
PO, with $\phi(\cdot,t)$ plotted on $\ga$; $\lam$ and $\phi$ over time $t$, 
for G and H 
showing the relaxation nature of the POs; and  $\kap$. For I we also 
plot the dynamic moving frame multipliers $\tau_x,\tau_y$, and $\rho$ (top right), 
and the space--time path of the center of mass of the cell (middle right). 
For G,H, $\tau_x,\tau_y$, and $\rho$ are zero ($<10^{-6}$ in modulus).
\label{al07h1}}}
\end{figure}

For the coupled system \reff{memo} we use $(\beta,\om,D,\del)=(1,1,0.1,0.3)$ 
and start with (fixed) $\al=0.7$, and continuation parameter $\zeta$. 
Figure \ref{al07f1} shows basic BDs, and samples of steady states; 
see Remark \ref{bdrem} for BD plot conventions, 
Fig.\ref{al07h1} for samples of POs, 
and \S\ref{ncsec} for numerical 
settings, i.e., mesh sizes and tolerances, and convergence 
experiments via mesh refinement in space and time.

 From \reff{zetc} we obtain that for increasing $\zeta$ the bifurcations 
from $\mS^1$ occur 
in the order $m=2$ (numerical $\zeta_{\text{BP}}=1.2723$, dark blue branch, 
samples A and B), $m=3$ (numerical $\zeta_{\text{BP}}=2.3125$, red 
branch, sample C), $m=4$ (numerical $\zeta_{\text{BP}}=3.7238$, 
violet branch, sample D), and $m=5$ (numerical $\zeta_{\text{BP}}=5.5266$, 
lilac branch, no sample), and so on.%
\footnote{\label{rbpfoot}\huc{The relative errors in $\zeta_{\text{BP}}$ 
between \reff{zetc} and the numerics in Fig.\ref{al07f1} 
are $<10^{-4}$ in all cases.}}
\footnote{\label{scfoot} The primary bifurcations from $\mS^1$ here 
are all subcritical for the given $\chi=50$, 
with a fold shortly after bifurcation. Here we are mostly 
interested in the behavior of the 
nontrivial branches away from $\mS^1$, but see 
\S\ref{aesec} for further discussion of the local 
 bifurcations from $\mS^1$ and their amplitude equations, and 
how that depends on $\chi$.} 
  The average morphogen coverage 
\huga{\label{phid}
\spr{\phi}:=\frac 1 {L_0}\int_\ga \phi(s)\dd s, 
\,\,\,\text{resp.}\,\,\, \spr{\phi}:=\frac 1 {L_0T}\int_0^T\int_\ga \phi(s,t)\dd s \dd t 
\text{ for POs with period $T$}, 
}
grows away from onset (a), while the surface tension Lagrange multiplier 
$\lam$ (resp. $\lam{=}\spr{\lam}{=}\frac 1 T \int_0^T \lam(t)\dd t$ for POs)  
decreases (b). In panel (c) we exemplarily show the translational Lagrange 
multiplier $\tau_y$, becoming non--zero in the drift bifurcation of the light 
blue branch.\footnote{The drift--in--$x$ multiplier $\tau_x$ is also nonzero on the light 
blue branch; no rotation is picked up in any of the (relative) steady 
states, i.e., $\rho=0$, but this will change on (some of) the PO branches, 
see Fig.\ref{al07h1}.} 

\brem\label{bdrem}{\rm 
In the BDs, thick/thin lines for steady state branches indicate 
linear stability/instability; for this, again the neutral translational 
and rotational modes are removed, as numerically they yield eigenvalues 
of order $10^{-6}$ with uncontrolled positive or negative real parts. 
For PO branches, no stability information 
is included in the BDs, and stability is discussed subsequently via DNS. 
BPs are indicated by $\circ$, HPs by $\diamond$, 
and sample points by $\bullet$, and labeled A, B, \ldots, respectively PO1, 
PO2,\ldots. The labeling of the BPs and 
HPs is Figure--wise, not branch-wise, i.e., HP3 in Fig.\ref{al07f1} 
is the third HP in the figure, but the first on the light blue branch.}
\eex\erem 

The $m{=}2$ branch is stable up to HP1 at $\zeta{\approx}4.2234$ where the 
PO branch with sample G  bifurcates.%
It is also 
again stable between HP2 at $\zeta\approx 5.3593$ (PO branch with H) 
and the BP at $\zeta\approx 5.8172$ (labeled BP2 in (b)), where the 
branch with sample E bifurcates, which connects 
to the $m=3$ branch at a BP at $\zeta\approx 3.4547$ (labeled BP1 in (b)), 
after which it continues to sample F and higher $\zeta$. Notably, 
the bifurcations at BP1 and BP2 are {\em drift bifurcations}, meaning that 
here the vesicles pick up a drift speed, indicated by the arrows in E and F, 
and by the translational Lagrange multiplier $\tau_y$ in (c). 
Naturally, BP1, where the $\CW_3$ 
branch and light blue branch meet, is a double BP: the kernel 
is spanned by vectors yielding the moving out or in of exactly 
one of the corners, yielding acute or obtuse triangles (like E), 
while distortion of the third corner corresponds to a linear 
combination of the two kernel vectors describing the other corners.%
\footnote{In the software, we just find a two dimensional kernel 
of this $D_3$ symmetry breaking bifurcation, and predictors for the 
three bifurcating branches must be found by (numerically) solving 
the so called algebraic bifurcation equations; see also 
\cite[\S2.5.3]{p2pbook} for a simple example of equivariant bifurcations 
with $D_3$ symmetry.} 
Note that the bifurcations at BP1 are transcritical, and at BP2 pitchforks, 
which again follows from symmetry. 
Finally, the $m=2$ branch itself 
connects to $m=4$ at $\zeta=5.4111$ (labeled BP3 in (b)), in a period 
(in $s$) halving pitchfork.

The $m{=}3$ branch is stable between BP1 and the fold at 
$\zeta{\approx}5.4304$. After this fold, it continues to 
larger $\spr{\phi}$, and for instance then features a HP at $\zeta\approx 4.049$ 
(not shown). For later comparison, here we remark that using DNS to test the 
linear stability indicated by the thick lines in Fig.\ref{al07f1} generally 
yields the following: for small perturbations of $\phi$ or $\ga$ 
from a stable steady state not close to a BP, FP or HP, we typically 
obtain convergence back to the perturbed steady state, but for larger 
perturbations the flow often goes to some ``run and tumble'' dynamics, 
discussed below in detail for $\al=3$. 
Here we summarize that already concerning just (relative) steady states 
the BD for \reff{memo} 
becomes rather complicated, with many secondary bifurcations 
from and interconnections of the primary branches. 

Figure \ref{al07h1} shows the three PO samples indicated in 
Fig.\ref{al07f1} as G,H,I. The top pictures 
show the shapes along the POs, with $\ga(s,t)$ colored by 
$\phi(s,t)$, where $s$ is the arclength along $\ga(\cdot,t)$; the red lines indicate the $s=0$ coordinate lines for the 
bottom $\kap(s,t)$ plots. The middle panels show the behavior 
of $\lam$ and $\spr{\phi}$ 
over $t$, indicating that both G and H are of relaxation oscillation type. 
The Lagrange multipliers $\tau_x,\tau_y$ and $\rho$ 
for the rotational and translational constraints are zero 
($<10^{-6}$ in modulus) for these POs, i.e., 
G and H are genuine breathers without hidden rotations or translations. 
The plots of $\kappa(s,t)$ at the bottom show how $\kap$ behaves along the POs; 
below we 
shall use such $\kap$ plots to illustrate more complicated dynamics. 
For I in (c), we additionally show $(\tau_x,\tau_y,\rho)$. Because the branch 
with I bifurcates 
from a branch of uniformly drifting relative steady states,  we already 
know that the averages $\spr{\tau_x}$ and $\spr{\tau_y}$ must be 
non--zero at bifurcation, and additionally the PO I picks up a 
time periodic non-zero rotational Lagrange multiplier $\rho$, 
but with $\spr{\rho}\approx 10^{-6}$. Thus, we obtain a net translation 
of $\ga$ over one period, but no net rotation; the last plot 
in (c) shows the space--time path of the center of mass (COM) of the 
moving vesicle, and its projection (red) into the $x$--$y$--plane.

All POs in Fig.\ref{al07h1} must be unstable. While the branches with G and H 
bifurcate at a stability loss/gain,  and hence could yield stable PO branches, 
these POs bifurcate subcritically in both cases. Now, running DNS 
(see \S\ref{timevo} for details) from 
small perturbations of any time slice $(\ga,\phi)_{t=t_0}$ 
of POs G or H, usually taking $t_0=0$, we get convergence back to 
the $\CW_2$ branch at the respective $\zeta$ value. 
The PO branch with I bifurcates 
from an already unstable steady state, and at this $\zeta$ 
range we do not find {\em any} stable steady states. Thus, we cannot 
expect convergence of perturbations of I to some steady state, 
and indeed such DNS yield self-intersection (and subsequent failure 
of the DNS) after rather short times. 

In summary, for $\al=0.7$, and altogether in the range $\al\in(0,1.5)$, say, 
and forcefully for $\al\le 0$, the model \reff{memo} only seems globally well-posed 
in general for small $\zeta$ ($\zeta<5$ for $\al=0.7$), where the 
dynamics is dominated by stable steady states with small $m$. Some of 
the steady state branches with small $m$ extend to larger $\zeta$, and 
higher $m$ branches bifurcate and exist at arbitrary large $\zeta$, 
but these are all unstable, and general ICs lead to self intersections 
and (numerical, and likely also analytical) blow--up. 
Moreover, we could not find any stable POs. 
However, both of these points change at larger $\al$, as shown in 
the next section.

\subsubsection{$\al=3$} \label{altw}
We now turn to the second regime, $\al=3$, for two reasons: (i) in this case, 
the order of primary bifurcations is 
$\zeta_3{\approx}3.1996
<\zeta_2{\approx}4.2010<\zeta_4{\approx}4.3829$  
\huc{(numerical values, with relative 
error $<10^{-4}$ compared to \reff{zetc})}, i.e., we have ``tuned'' 
the primary 
bifurcations to the $\CW_3$ branch bifurcating first. (ii) this 
naturally has further consequences for the secondary bifurcations, 
and inter alia we now find stable $\Z_2\times\Z_2$ symmetric breathers  
bifurcating from the $\CW_4$ branch. 

\begin{figure}[ht]
\centering
\btab{l}{{\sm (a)}\hs{90mm}{\sm (b)} \\
 \hs{-5mm}\ig[width=0.58\tew]{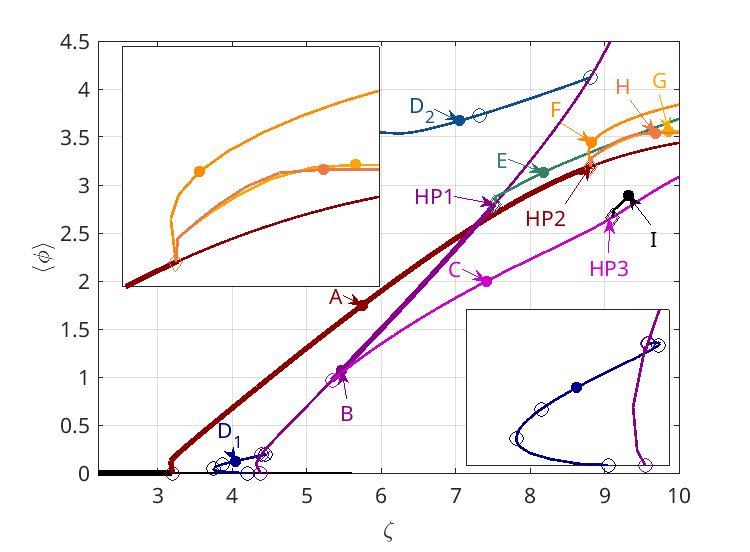}
\hs{-5mm}\ig[width=0.35\tew]{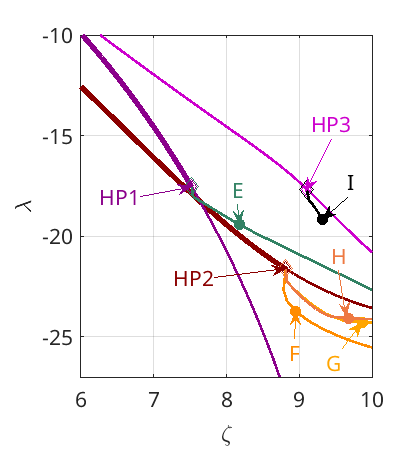}\\[-4mm]
{\sm (c)}\\
\hs{0mm}\ig[width=0.18\tew]{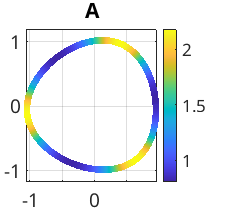}\hs{0mm}\ig[width=0.19\tew]{3/B}
\hs{0mm}\ig[width=0.21\tew]{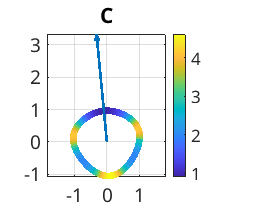}\hs{0mm}\ig[width=0.17\tew]{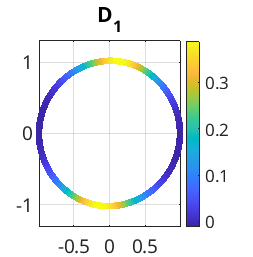}
\hs{0mm}\ig[width=0.17\tew]{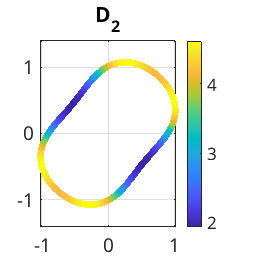}
}\\
{\sm 
\begin{tabular}{l|cccccccccc} \hline 
&A&B&C&D$_1$&D$_2$&E&F&G&H\\
$\zeta$&4.8085&5.4703&7.4205&4.0439&7.0516&8.1810&8.8317&9.8624&9.6780\\
\hline
\end{tabular}}
\vs{0mm}
\caption{\small{Basic BD  for \reff{memo} 
with parameters from \reff{bps} and $\al=3$, and 
steady state samples (C a relative steady state, as indicated 
by the arrow); see Fig.\ref{al3c2} for PO samples. }  
\label{al3c1}}
\end{figure}

Figure \ref{al3c1} shows a basic BD of the three primary branches 
$\CW_3$ (brown, 
with sample A), $\CW_2$ (blue, with sample D$_1$, see also the 
zoom inset at the bottom right),  
and $\CW_4$ (violet, with sample B, 
also already previewed in Fig.\ref{f2a}(a)), two secondary steady 
state branches 
(with samples C and D$_2$), 
and five PO branches (green, orange, and black, with labels E--I, see also zoom inset 
at the top left). 
We start with the $\CW_2$ branch, 
which here bifurcates more strongly subcritical 
and on the other hand after two folds reconnects with the $\CW_4$ branch, 
see zoom inset. The $\CW_2$ branch is hence not of major interest to us,  
and just illustrates the strong effect the coupling with $\phi$ can 
have also close to $\mS^1$. For completeness, we also note 
that there is a second $\CW_2$ branch bifurcating from $\CW_4$ at large amplitude, 
with sample D$_2$, which however is unstable throughout.

The $\CW_4$ branch 
becomes stable in a double BP at $\zeta\approx 5.3439$.  Here 
the kernel is spanned by vectors describing the moving out and in 
of two opposite ``corners'', which can occur in two orientations. 
The magenta  
branch with sample C is one (direction of one) of the bifurcating branches, 
namely with the top corner moving in (and the bottom slightly moving out), 
such that the $D_4$ symmetry is broken 
to $\Z_2$. As a result, this is again a drift bifurcation, with the picked 
up speed indicated by the arrow in sample C. On this branch there 
is a Hopf point, and the POs bifurcating there also pick up some rotation, 
similar to sample H discussed below.  

The $\CW_4$ branch then loses stability again at 
$\zeta\approx 7.5216$ 
in a {\em simple} HP labeled HP1. From symmetry it then follows that 
a single branch of 
POs must bifurcate, on which solutions alternate between two $\Z_2\times\Z_2$ 
symmetric shapes, i.e., a breather, without drift or rotation, see 
sample E (already previewed in Fig.\ref{f2a}(b)). 
As this Hopf bifurcation is supercritical, we expect these breathers 
to be stable, at least near onset, although this stability is difficult 
to analyze due to the always present zero modes from spatial translations and 
rotations. Numerically, we find the bifurcating solutions 
to be {\em metastable} under DNS in the labframe, 
i.e., only stable over long transients, 
but stable under constrained DNS, 
see Fig.\ref{al3f}, and the discussion there. 

The situation is more difficult at the {\em double HP} on $\CW_3$ labeled 
HP2 at $\zeta{\approx}8.8149$, where $\CW_3$ loses stability. The Hopf eigenspace 
here is spanned by two vectors each describing the dominant motion of  
just one corner, and by symmetry a dominant motion of the third corner 
can be described by a linear combination of the first two. From $D_3$ 
symmetry it follows 
that at least (and generically exactly) 
three PO branches must bifurcate, modulo conjugacy, see Remark \ref{d3sym}.

\begin{figure}[H]
\btab{l}{{\sm (a)}\\
\hs{-0mm}\ig[width=0.22\tew]{3/E}
\rb{12mm}{\btab{l}{\ig[width=0.22\tew]{3/El}\\[-4mm]
\ig[width=0.22\tew]{3/Epa}}}
\ig[width=0.26\tew]{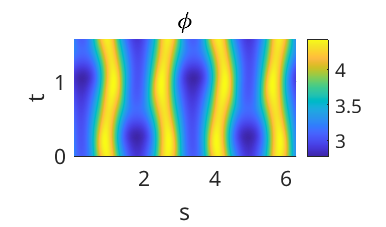}\ig[width=0.26\tew]{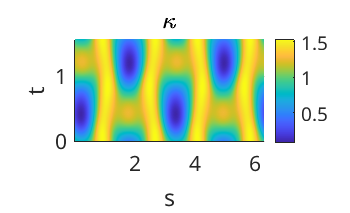}\\
\hline
\btab{l|l}{{\sm (b)}&{\sm (c)} \\
 \hs{-0mm}\btab{l}{
\ig[width=0.22\tew]{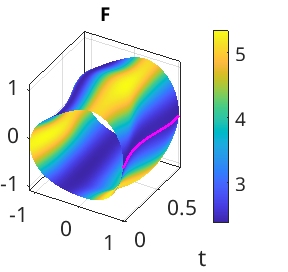}\\
\ig[width=0.22\tew]{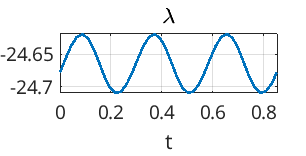}\\
\ig[width=0.22\tew]{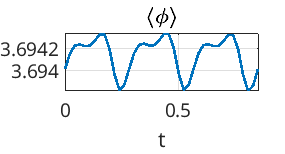}}
\hs{-5mm}\btab{l}{
\ig[width=0.24\tew]{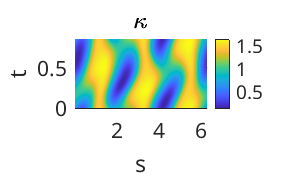}\\
\ig[width=0.21\tew]{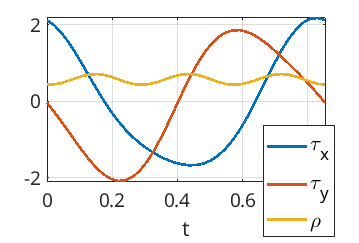}\\
\ig[width=0.21\tew]{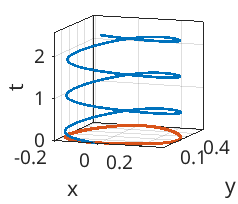}}
&
\hs{-0mm}\btab{l}{
\ig[width=0.22\tew]{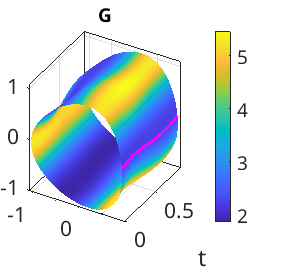}\\
\ig[width=0.21\tew]{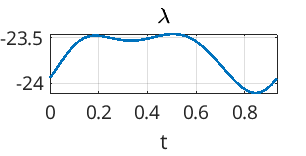}\\
\ig[width=0.21\tew]{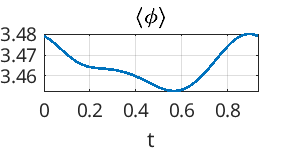}}
\hs{-5mm}\btab{l}{
\ig[width=0.24\tew]{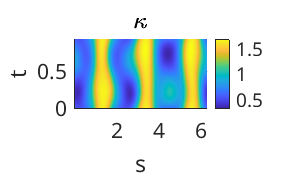}\\
\ig[width=0.21\tew]{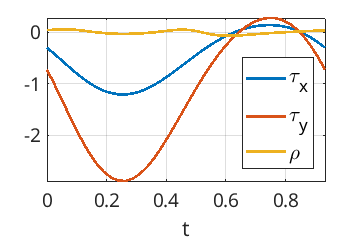}\\
\ig[width=0.17\tew]{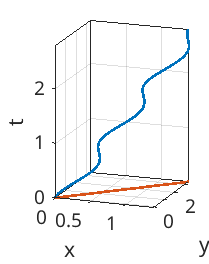}}
}\\
\hline
\rb{25mm}{\hs{-2mm}{\sm (d)}}
\hs{-0mm}\ig[width=0.24\tew]{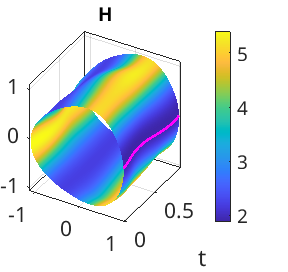}
\rb{14mm}{\hs{-2mm}\btab{l}{\ig[width=0.24\tew]{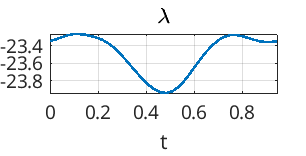}\\[-4mm]
\ig[width=0.24\tew]{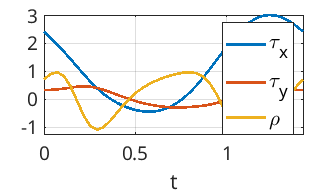}}}
\hs{-2mm}\ig[width=0.27\tew]{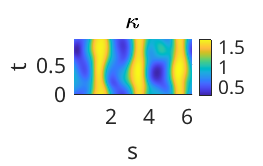}
\ig[width=0.24\tew]{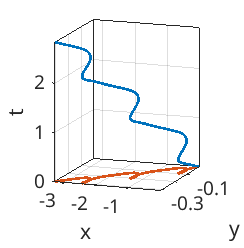}
}
\vs{-3mm}
\caption{\small{PO samples for Fig.\ref{al3c1}; (a) $\Z_2\times\Z_2$ symmetric 
breather; (b-d) (relative) POs bifurcating from the $D_3$ symmetric 
HP2. See text for details. }  
\label{al3c2}}
\end{figure}

Figure \ref{al3c2} illustrates four PO samples marked E--H in Fig.\ref{al3c1}. 
As already said, E in (a) is a ``breather'', and in particular the 
translational and rotational 
multipliers \amc{$\tau_x,\tau_y$ and $\rho$ }
are ${<}10^{-6}$ (in modulus) and hence not shown. 
The two $\Z_2$ 
symmetries (along $s{=}0$ and $s{=}\pi$), together with 
the $(s,t)\mapsto (s{+}\pi/4,t{+}T/2)$ symmetry can clearly be seen. 
For PO F in (b), \amc{$\tau_x,\tau_y,\rho{\ne}0$}, but still $\spr{\tau_x}{=}0$ 
and $\spr{\tau_y}{=}0$, and hence there is no net motion over one period $T$; 
the space--time plot of the COM in the last plot in (b) forms a helix, 
and the projection in the $x$--$y$ plane a closed loop. We believe that 
the best name for this kind of PO is ``modulated traveling wave'', 
see also the $\kap$ plot in (b). 
In contrast, for PO G in (c) we see a ``stick--slip'' motion of the 
cell in the lab frame. Here, the top ``corner'' at $s{=}s_0{\approx}1$ oscillates 
most strongly, and the two other corners follow with equal phase, and hence 
$\ga$ alternates between an obtuse and an acute triangle. 
Letting $\vec{v}=(\cos(s_0),\sin(s_0)){\approx}(1,2)/\sqrt{5}$ 
this yields motion in direction $\vec{v}$, sticking (in fact 
slightly reversing) when $\ga$ is obtuse, and slipping when $\ga$ is acute. 
See the COM mass plot, with projection in the $x$--$y$ plane a line in direction 
$\vec{v}$. In (d), the third PO H bifurcating from the \amc{double} HP HP2 can be seen 
as a superposition of F and G. Moreover, we note that the PO I bifurcating 
from the branch of rigidly translating vesicles (see sample C) is 
very similar to the PO H.

\brem\label{d3sym}{\rm 
The three branches with samples F,G and H  bifurcating at HP2 in Fig.\ref{al3c1} 
and illustrated in Fig.\ref{al3c2}(b--d) 
represent the three conjugacy classes 
(equivalency of POs under rotation by $2\pi/3$) of branches 
bifurcating at double HPs in the $D_3$ symmetric case, in the 
one-dimensional fixed point subspaces of the $D_3$ group action.  
In the numerics, to search for different branches of POs 
bifurcating at HPs of higher multiplicity $n\ge 2$, we optionally pass a 
vector $(z_1,\ldots,z_n)$ of coefficients for the different 
eigenvectors $\Psi_1,\ldots,\Psi_n$ to the pertinent frequency $\om_H$ 
to the branch switching routine%
\footnote{see also \cite[\S7.2.3]{p2pbook} for an example 
of double HPs in an example with $D_4$ symmetry}, which then 
creates the predictor 
\huga{
U_{\text{pred}}(s,t;z)=U_{H}(s)+\ddsi[(z_1\Psi_1(s)+\ldots
+z_n\Psi_n(s))\er^{\ri\om_H t}+\cc], 
}
for branch switching to a PO branch. In the double $D_3$ HP cases, 
different choices of $z_1,z_2$ then always takes us to one of the branches with 
F,G or H, modulo conjugacy, i.e., modulo rotation of $\ga$ (and hence the 
COM paths) by $2\pi/3$, or, in rare cases, non--convergence of the 
corrector for the PO computation.   This agrees with the general 
theory for $D_3$ equivariant Hopf bifurcations,  see for instance 
\cite{DP05} and the references therein, which yields that 
there generically bifurcate exactly three (not just at least three) 
branches (modulo conjugacy) of POs at double HPs. }
\eex\erem

\subsubsection{Stability and DNS}
\begin{figure}[H]
\btab{llll}{
{\sm (a)}&{\sm (b)}&{\sm (c)}&{\sm (d)}\\[-0mm]
\hs{-4mm}\rb{20mm}{\btab{l}{\hs{0mm}\ig[width=0.17\tew]{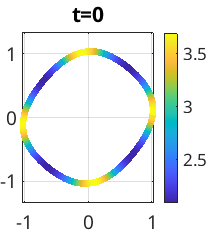}
\hs{-4mm}\ig[width=0.17\tew]{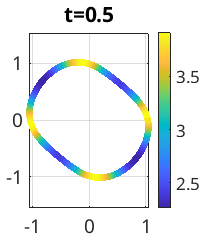}\\
\ig[width=0.17\tew]{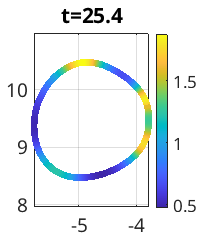}
\hs{-4mm}\ig[width=0.17\tew]{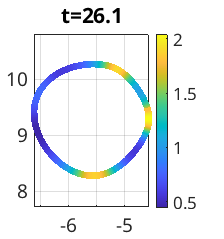}}}
&\hs{-5mm}\rb{18mm}{\btab{l}{\ig[width=0.24\tew,height=30mm]{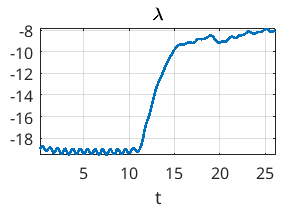}\\
\ig[width=0.24\tew]{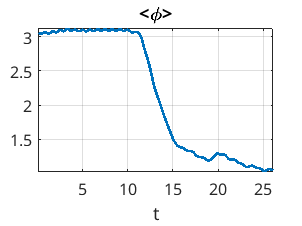}}}
&\hs{-6mm}\ig[width=0.22\tew]{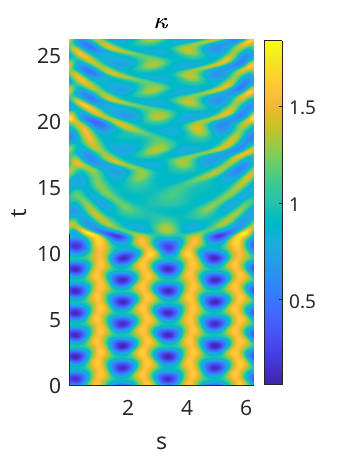}&\hs{-3mm}\ig[width=0.2\tew]{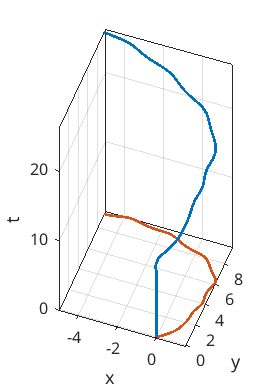}
}
\\[-3mm]
\btab{l}{
{\sm (e)}\\
\hs{-0mm}\ig[width=0.26\tew]{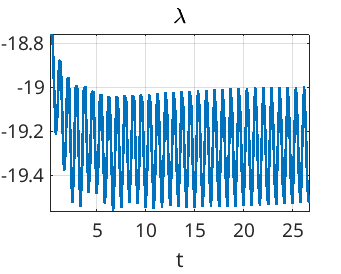}
\hs{-0mm}\ig[width=0.26\tew]{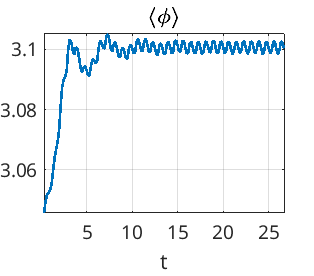}
\hs{-1mm}\ig[width=0.27\tew]{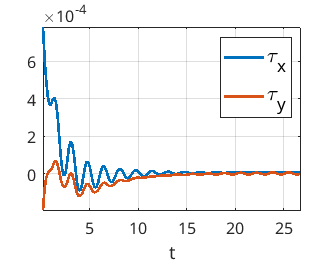}
\ig[width=0.22\tew]{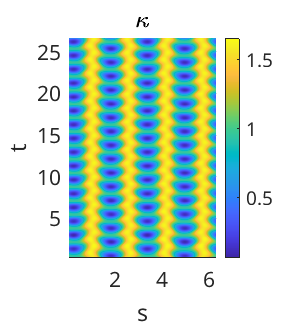}
}

\vs{-4mm}
\caption{{\sm (a) Snapshots from  
DNS from perturbations of the PO E time $t=0$ slice, 
$X|_{t=0}=X_0+a\cos(3\vt)N$ 
with amplitude $a=0.01$. (b) time--series of $\lam$ and $\spr{\phi}$. 
(c) (together with (a)) illustrates the convergence to a tumbling vesicle 
with a source--sink pair 
(trigger wave) for $\kap$. (d) shows the COM path of the vesicle.
(e)  cDNS. All results with {\tt implicitEuler}, but same with {\tt RADAU} (in (a-d) 
up to $t=20$, then failure). }
\label{al3f}}
\end{figure}

In Fig.\ref{al3f}(a-d) we present 
DNS {\em in the lab frame}, i.e., using \reff{dae1}, with a 
perturbation of the $t=0$ time slice of PO E as IC. The PO E 
in principle should be stable as the branch bifurcates supercritically 
at the stability loss of the $\CW_4$ branch (again after removing 
the neutral modes, see Remark \ref{bdrem}). However, in our 
straightforward 
DNS, the solutions only appear to 
be metastable, i.e., stable on long but 
finite timescales (up $t=t_0=14$, say). We believe 
this is due to (unavoidable) numerical errors 
interacting with the neutral translational and rotational modes. 
  Once these get sufficiently 
strongly excited (after $t=t_0$), a net motion sets in (see the 
COM path in (d), and the axis ticks in (a)), and in its comoving 
frame the solution goes to a roughly periodic (but with 
decreasing amplitude, see (b) after $t=25$) source--sink pair dynamics for 
$\kap$ (and $\phi$) on $\ga$, somewhat similar to the trigger waves 
from Fig.\ref{ef1}J--L. 
For smaller amplitude initial perturbations, and/or for 
POs closer to HP1, i.e., for smaller amplitude non--moving breathers, 
we obtain longer metastability, but eventually we always end up 
with source--sink pairs, in which moreover 
$\lam(t)$ always goes to a higher level, and 
$\spr{\phi(t)}$ to a lower level.

In the constrained DNS (cDNS, cf.\,\reff{dae2}) 
in Fig.\ref{al3f}(e), trading spatial translations for 
translational Lagrange multipliers, the PO E becomes genuinely stable. 
This shows some inconsistency between \reff{dae1} and \reff{dae2}, 
which however only appears after many periods of oscillations (also 
depending on the initial perturbation), and as already said we believe 
that over long times the unconstrained DNS in (a-d) might be less reliable than 
the cDNS in (e) due to 
uncontrolled neutral translational (and rotational) modes in (a-d).%
\footnote{In cDNS, here we only use the translational 
multipliers $\tau_x$ and $\tau_y$, and no constraint for rotations; 
however, these do not get excited, and we obtain the same results for 
cDNS including such rotational constraints. 
Moreover, we checked that the cDNS here and all other cases 
gives the same results when using a simple implicit Euler scheme 
(which can be overdamping) for the DAEs, and when using the high 
order DAE--suitable code {\tt RADAU}, see the discussion after \reff{dae2}. 
On the other hand, {\tt RADAU} typically only works for rather short time 
($t\le 2$, say) 
in DNS (without translational and rotational constraints), and then 
yields excessive step--size reductions and eventually non--convergence.  
Therefore, the presented (lab frame) DNS are always done with implicit Euler.}

\begin{figure}[h]
\btab{l}{
{\sm (a)}
\\[-6mm]
\ig[width=0.22\tew]{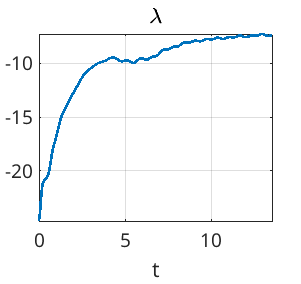}\hs{-0mm}\ig[width=0.22\tew]{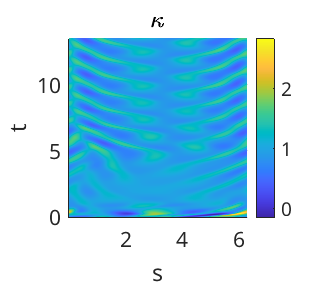}
\rb{2mm}{\ig[width=0.21\tew]{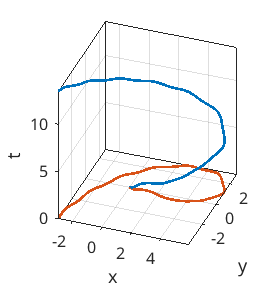}}\hs{4mm}\ig[width=0.22\tew]{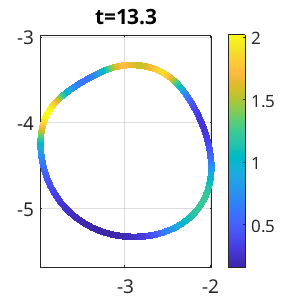}\\
{\sm (b)}
\\[-4mm] 
\ig[width=0.22\tew]{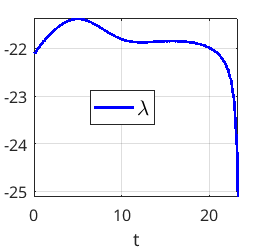}\ig[width=0.35\tew]{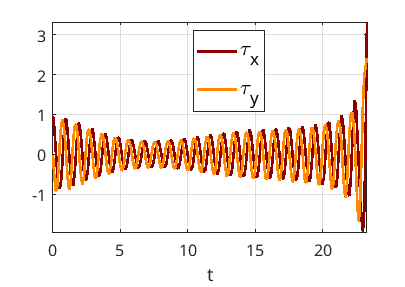}
\ig[width=0.24\tew]{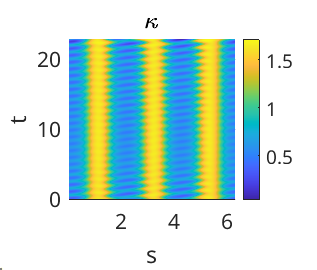}\ig[width=0.22\tew]{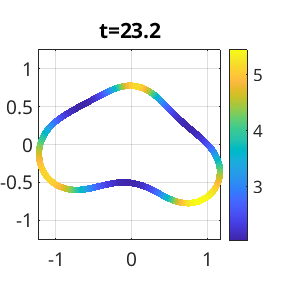}\\
{\sm (c)}
\\[-4mm]
\ig[width=0.24\tew]{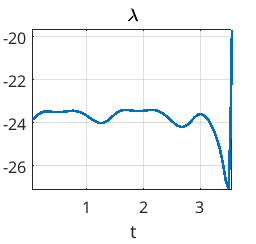}\hs{-0mm}\ig[width=0.27\tew]{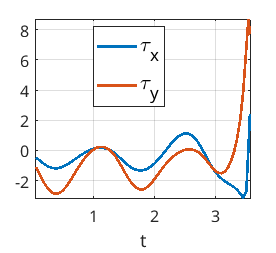}\hs{-0mm}\ig[width=0.24\tew]{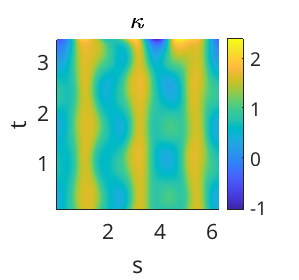}\rb{5mm}{\ig[width=0.23\tew]{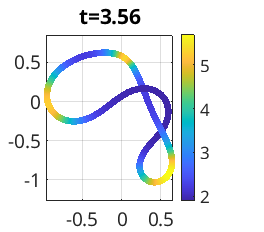}}
}

\vs{-4mm}
\caption{{\sm (a) DNS from F, emergence of traveling vesicle with trigger waves.
(b) cDNS from F, metastable behavior, but eventual self--intersection 
in implicit Euler, and failure of RADAU.
(c) cDNS from G,  unstable, self intersection (and failure of RADAU) at small 
$t$. 
}
\label{al3f3}}
\end{figure}

For initial perturbations of the $t=0$ slice 
of the  $D_3$ symmetric POs F,G, and H we obtain the following, 
see Fig.\ref{al3f3}. 
Under DNS (a), F is unstable and after a rather short transient 
yields trigger waves. On the other hand, for cDNS, F shows 
some metastability: starting from a moderate initial perturbation 
(of amplitude 0.1 times the amplitude of the PO F), the solution 
shows periodic behavior up to $t\approx 20$, which is best seen 
in the behavior of $\tau_x$ and $\tau_y$, and the $s$--$t$ plot of $\kap$. 
Subsequently, the shapes 
of the oscillations deteriorate, and at $t\approx 23.2$ {\tt RADAU} fails, 
while implicit Euler leads to self--intersection and subsequent failure. 
Finally, DNS from G,H are genuinely unstable in both, DNS and cDNS. 
In DNS, they again yield convergence to trigger waves, and 
in cDNS self intersection (c, for G), and similar for H. 

\section{Experiments in 2D}\label{3dsec}
Having obtained some understanding of the multitude of possible steady states, 
POs, and dynamics, for \reff{memo} in 1D in \S\ref{1Dsec}, we come back 
to the 2D problem \reff{mm1}, with one sample PO  
already given in Fig.\ref{f2a}(c). 
The analytical and/or numerical 
classification and stability 
problem of possible steady states and POs in 2D in dependence of 
parameters becomes {\em much} harder, mostly due to the 
higher multiplicity of BPs and HPs in 2D, also related to the larger 
symmetry group $O(3)$ of rigid body motions of the sphere compared to 
$O(2)$ for the circle. 
Thus, here we mainly restrict to one basic BD to which the PO from Fig.\ref{f2a}(c) 
belongs, and give some DNS experiments. 

\begin{figure}[H]
\btab{l}{
\btab{ll}{{\sm (a)}&{\sm (b)}\\
\rb{0mm}{\btab{l}{\hs{-0mm}\ig[width=0.4\tew]{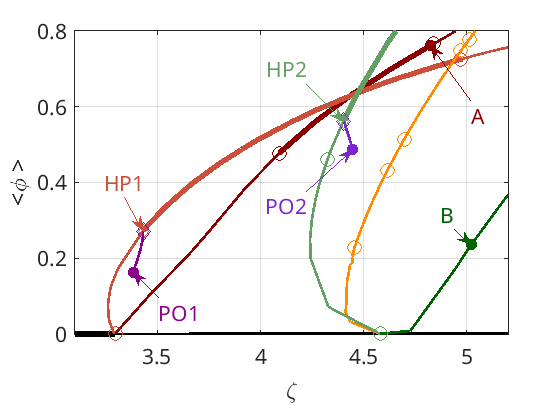}\\[-7mm]
\ig[width=0.4\tew]{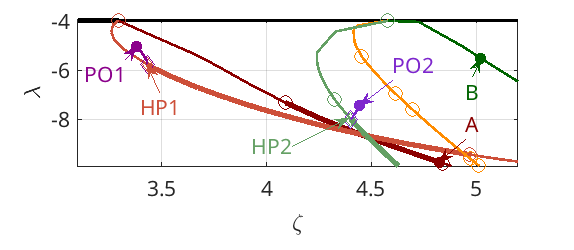}}}&
\hs{0mm}\rb{0mm}{\btab{ll}{\ig[width=0.23\tew]{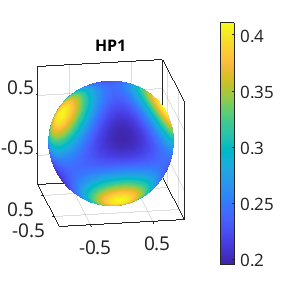}
\hs{2mm}\ig[width=0.23\tew]{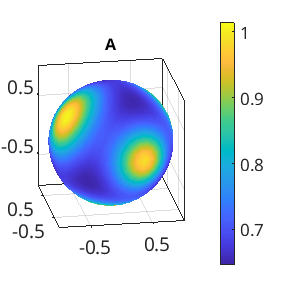}\\
\ig[width=0.22\tew]{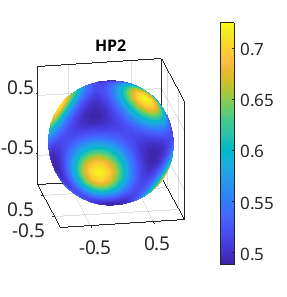}\hs{2mm}\ig[width=0.22\tew]{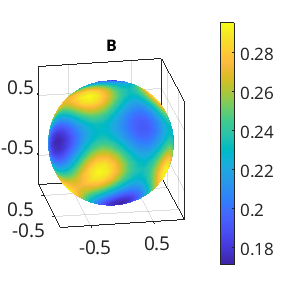}
}}
}\\[-4mm]
{\sm (c)}\\[-5mm]
\hs{10mm}\ig[width=0.22\tew,height=33mm]{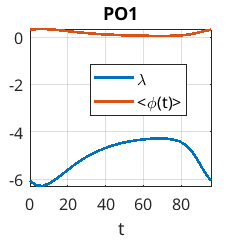}\hs{5mm}\ig[width=0.2\tew]{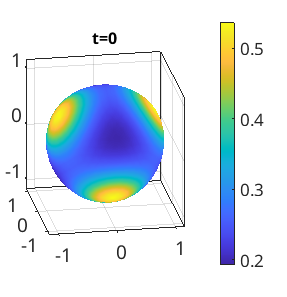}
\ig[width=0.2\tew]{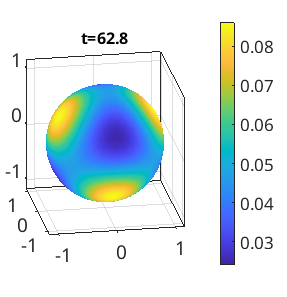}\\
{\sm (d)}\\[-8mm]
\hs{10mm}\ig[width=0.22\tew]{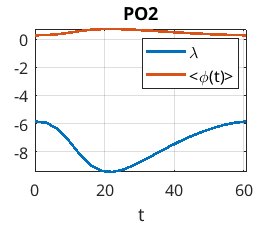}\hs{5mm}
\ig[width=0.2\tew]{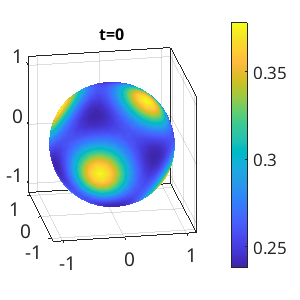}\ig[width=0.2\tew]{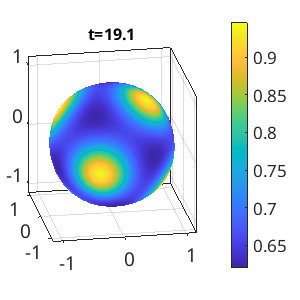}
}\\
\centering
{\sm 
\begin{tabular}{l|cccccccccc}   \hline 
&HP1&HP2&A&B&PO1&PO2\\
$\zeta$&3.4302&4.4020&4.8235&5.0219&3.3815&4.4458\\ \hline
\end{tabular}}
\vs{-0mm}
\caption{{\sm  BDs (a) and steady state samples (b)  
for \reff{mm1} with $(\beta,\om,D,\del){=}(1,1,0.1,0.05)$ and $\al=2$. 
First bifurcation from the sphere at $\zeta{\approx}3.29$ to $l{=}3$ spherical 
harmonics (7 dimensional kernel), 
and second bifurcation at $\zeta{\approx}4.57$ to 
$l{=}4$ (9 dimensional kernel). At both BPs we only show selected 
bifurcating branches. (c,d) PO1 and PO2 samples. 
\label{3f5a}}}
\end{figure}

Panel (a) in Fig.\ref{3f5a} shows selected steady state branches for 
$\al=2$, and two PO branches. Steady state samples are given in (b); 
(c) shows a PO sample from the first PO branch with tetrahedral 
symmetry, already previewed in Fig.\ref{f2a}, and (d) shows a sample 
from PO2 with hexahedral symmetry. 
The bifurcations from the 
sphere happen in spherical harmonic order $l{=}3$ at $\zeta{\approx}3.2945$ (BP1), 
$l{=}4$ at $\zeta\approx{4.5732}$ (BP2) and $l{=}5$ at 
$\zeta\approx 6.3433$ (not shown). 
For $l{=}3$ at BP1, the kernel is hence 7-dimensional, and 
for $l{=}4$ at BP2 it is 9-dimensional, 
and at BP1 we only compute two bifurcating steady state branches, 
and three at BP2. 
We start the discussion with the brown regular tetrahedron branch with HP1. 
This branch becomes stable at the Hopf bifurcation of breathing 
tetrahedra, and loses stability in a steady state bifurcation 
at $\zeta{\approx}4.97$ to ``distorted'' tetrahedra (not shown). 
The branch with A, with solutions with 5 unequal spots 
bifurcates simultaneously with the tetrahedral branch at BP1, 
together with several other branches. Similarly, we show 
three steady state branches bifurcating at BP2, of which the 
green one is most interesting as it gains stability at HP2, 
where the PO2 branch of breathing hexahedra bifurcates.
Both, HP1 and HP2 are simple, and in summary, 
the bifurcating PO branches of breathing tetrahedra and hexahedra 
are analogous to the breathing $D_4$ branch with sample E from 
Fig.\ref{al3c2}(a). In particular, the three translational 
 and the two rotational  Lagrange multipliers are all zero on 
the branches with PO1 and PO2.

As in 1D, the (linearized) 
stability of steady states is indicated by thicker lines, and we checked 
that this fully agrees with stability obtained from DNS starting 
with small perturbations of linearly stable steady states. 
For the POs, we may expect stability of PO1 
as the branch bifurcates supercritically from the tetrahedral branch, 
while PO2 must be unstable as its branch bifurcates subcritically. 
To assess dynamics near these  POs, and more general dynamics also 
yielding further stable steady states at larger $\zeta$, 
we again resort to DNS, similarly as in 1D 
as DNS in the lab frame, and as cDNS, now 
with additional translational Lagrange multipliers $(\tau_x,\tau_y,\tau_z)$.%
\footnote{Again, the rotational neutral modes appear to be less dangerous here, 
and hence we do not add constraints and Lagrange multipliers for these in 
cDNS.}

\begin{figure}[H]
\btab{l}{
{\sm (a)}\\[-6mm]
\ig[width=0.24\tew]{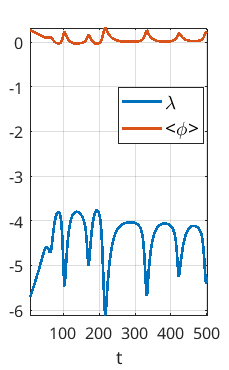}\hs{-3mm}
\ig[width=0.24\tew]{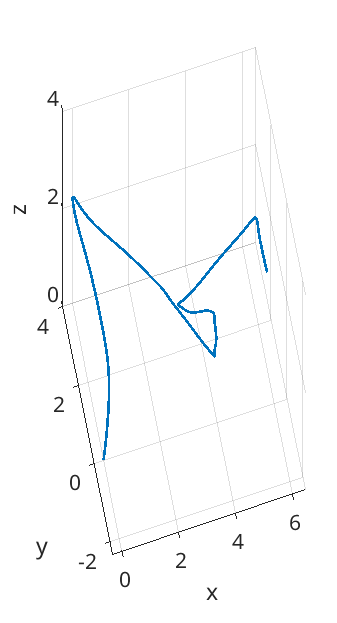}\hs{-3mm}
\hs{4mm}\rb{36mm}{\btab{l}{
\ig[width=0.25\tew]{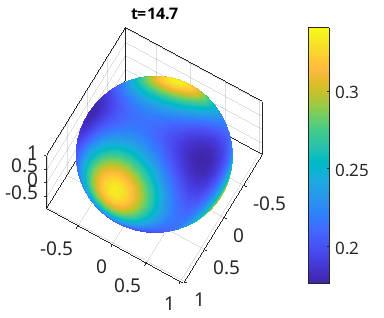}\hs{-3mm}
\ig[width=0.25\tew]{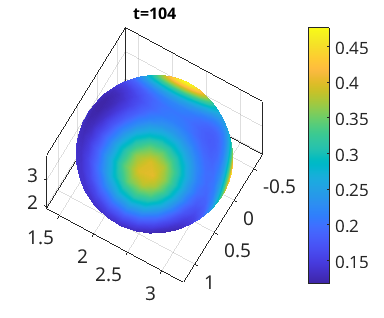}\\
\ig[width=0.25\tew]{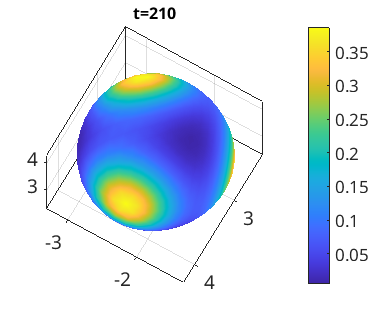}\hs{-3mm}
\ig[width=0.25\tew]{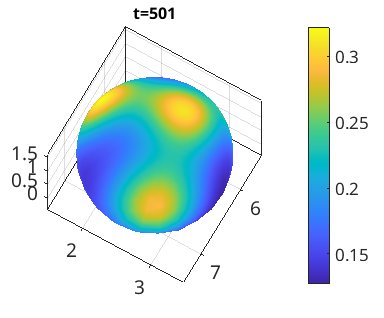}}}\\
{\sm (b)}\\[-2mm]
\ig[width=0.24\tew]{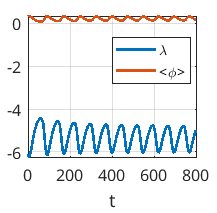}\hs{-0mm}
\ig[width=0.24\tew]{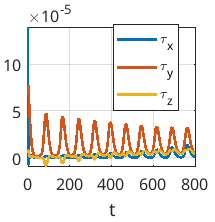}\hs{-1mm}
\ig[width=0.24\tew]{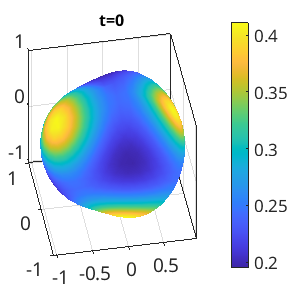}\hs{-1mm}
\ig[width=0.24\tew]{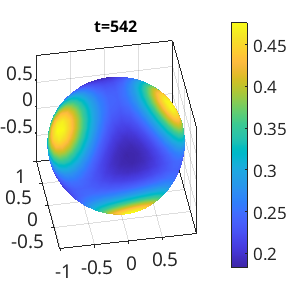}
}

\vs{-2mm}
\caption{{\sm  (a) DNS starting near HP1 from Fig.\ref{3f5a}; 
drifting and tumbling vesicle, with numbers of $\phi$ spots alternating between 
4 and 5. 2nd plot is the COM path. 
(b) cDNS starting from a (rather large) perturbation 
of the $t=0$ slice of PO1, showing stability of PO1 under cDNS.   
\label{3f5b}}}
\end{figure}

In Fig.\ref{3f5b} we contrast the  labframe DNS in (a) with cDNS in (b). 
The time series for $\lam$ and $\spr{\phi}$ in (a) show some 
similarity (in period and amplitude) to those of PO1, but like in 1D we 
also get some irregular tumbling and motion of the vesicle 
(COM path in 2nd plot), 
together with shape deformations which lead to the 
appearance, shifting and merging (e.g., near $t=500$) of the 
spots of $\phi$ in the snapshots on the right. However, 
we mostly observe 4 spots in $\phi$, and roughly tetrahedral symmetry. 
In contrast, in (b) we take a rather large perturbation of the $t=0$ time slice 
from PO1, and get convergence back to PO1, with the 
translational Lagrange multipliers $(\tau_x,\tau_y,\tau_z)$ staying 
small for very long times (up to $t=2000$). 

\begin{figure}[ht]
\iffalse
\btab{l}{{\sm (a)}\\[0mm]
\ig[width=0.23\tew]{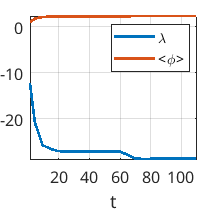}\hs{-0mm}
\rb{3mm}{\ig[width=0.26\tew]{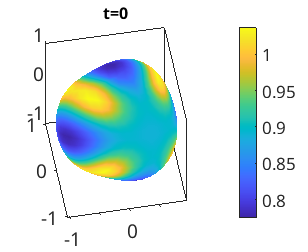}\hs{-1mm}
\ig[width=0.26\tew]{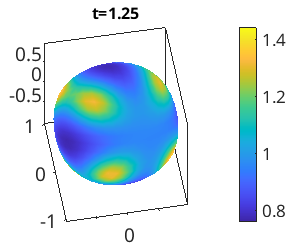}\hs{-1mm}
\ig[width=0.26\tew]{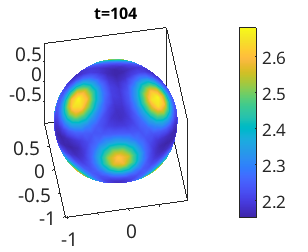}}\\
{\sm (b)}\\[-0mm]
\ig[width=0.23\tew]{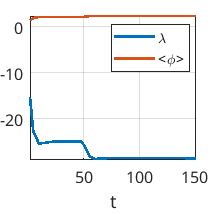}\hs{-0mm}
\rb{3mm}{\ig[width=0.26\tew]{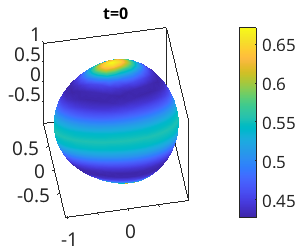}\hs{-1mm}
\ig[width=0.26\tew]{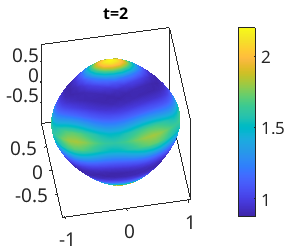}\hs{-1mm}
\ig[width=0.26\tew]{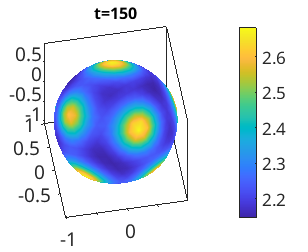}}}
\else
\btab{l}{{\sm (a)}\\[-3mm]
\hs{-3mm}\ig[width=0.2\tew]{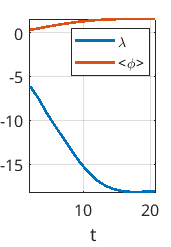}\hs{-0mm}
\rb{0mm}{\ig[width=0.2\tew]{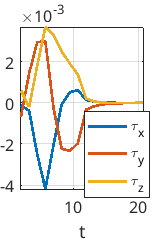}\hs{-1mm}
\ig[width=0.33\tew]{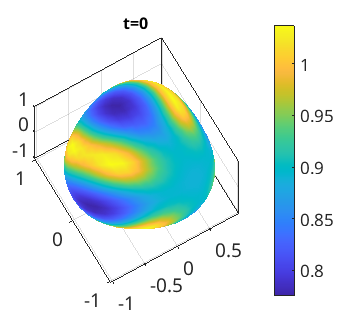}\hs{-1mm}
\ig[width=0.33\tew]{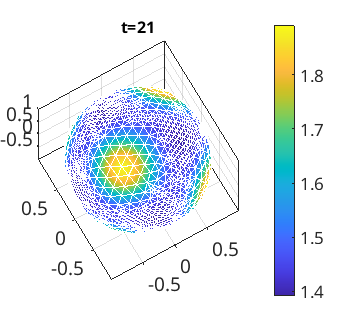}}\\
{\sm (b)}\\[-3mm]
\hs{-3mm}\ig[width=0.2\tew]{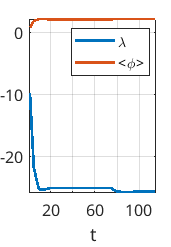}\hs{-0mm}
\rb{0mm}{\ig[width=0.2\tew]{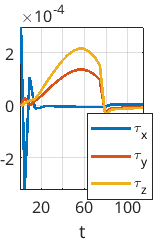}\hs{-1mm}
\ig[width=0.33\tew]{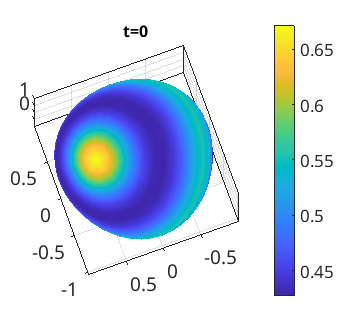}\hs{-1mm}
\ig[width=0.33\tew]{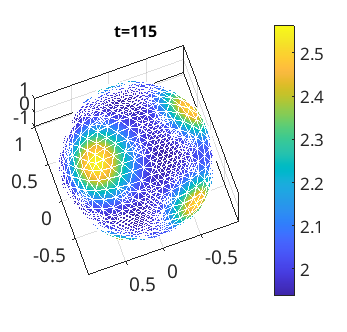}}}
\fi
\vs{-2mm}
\caption{{\sm  (a) cDNS from perturbation of B ($\zeta=5.0219$) from 
Fig.\ref{3f5a} but at $\zeta{=}6$; convergence to octahedron 
(two square pyramids, 6 spots in $\phi$). 
(b) cDNS from perturbation of a solution in Fig.\ref{3f5a} on the 
orange branch at 
$\zeta{\approx}4.4458$ (see $t{=}0$ plot in (b)), 
but at $\zeta{=}7$; convergence to decahedron (two pentagonal pyramids, 
7 spots in $\phi$). Last plots show refined meshes. 
\label{3f6}}}
\end{figure}

Finally, in Fig.\ref{3f6} we use cDNS at larger $\zeta$, 
and obtain convergence to stable steady states (of higher 
spherical harmonics type), and this seems to be the generic behavior 
for the given parameters, up to $\zeta=12$. However, the associated steady 
state branches are rather difficult to continue at larger $\zeta$, 
due to (faulty) loss of symmetry during the continuation, indicating branch 
jumping.  Moreover, the translational multipliers 
$\tau_x,\tau_y$, and $\tau_z$, 
generally converge to $\CO(10^{-4})$ values, reflecting some 
asymmetry in the obtained meshes, but this can be reduced via 
adaptive mesh--refinement; in (a) we do this after the 
solution has roughly reached the octahedral shape at $t=10$, and 
in (b) at $t=75$ after reaching the decahedral shape, in both cases yielding 
$\tau_x,\tau_y,\tau_z=\CO(10^{-5})$ at the respective ends of the DNS.%
\footnote{\label{mrfoot} \huc{Here we select elements with 
{\em small} $\phi$ for refinement, which turned out to give a stronger 
reduction in $\tau_x,\tau_y,\tau_z$ than, e.g., the converse, i.e., 
refinement at large $\phi$, and which also turned out to be more 
robust for using the final states in Fig.\,\ref{3f6} as initial 
guesses for continuation as steady states to larger $\zeta$, see 
\S\ref{ncsec}. In any case, the mesh adaption acts to somewhat 
artificially (but mildly) enforce the symmetry class. 
See also, e.g., \cite{ZHL17} for numerical methods to enforce a given symmetry 
on a solution branch in the related problem of two--phase vesicles. }}

\section{Discussion and outlook}\label{dsec}
We studied what we 
believe to be one of the simplest geometric mechanochemichal models, 
which couples the Helfrich gradient flow for the vesicle membrane 
$X$ (resp.$\ga$ in 1D) with a reaction--diffusion equation 
for a single morphogen $\phi$ on $X$ (resp.~on $\ga$), via 
the dependence of the spontaneous curvature as $c_0=\al+\beta\phi$, 
and conversely via the kinetics $\ds \zeta \frac{f(\kap)}{\om+f(\kap)}$ 
for $\phi$. 
We focused on the 1D case, which is already very rich in 
steady and dynamical pattern formation, and for 2D 
restricted to an outlook. 

Our primary continuation parameter was $\zeta$, 
and  guided by the linear stability 
analysis at the circle in \S\ref{locstab} we took two different 
values for $\al$ with remaining parameters fixed: 
For small $\al=0.7$ (\S\ref{alon}), the primary loss of stability 
of $\mS^1$ is 
wrt to spatial wave number $m=2$ steady patterns, then $m=3, m=4, \ldots$. 
However, larger $\al$ can be used to sort larger $m$ to the front, 
and in \S\ref{altw} we used $\al=3$ yielding the order $m=3,2,4,\ldots$.   
In both cases, we find HPs on the primary bifurcating branches, some 
of these of higher multiplicities, and 
from these HPs bifurcations of POs, which essentially come 
in two variants: {\em breathing} vesicles, i.e., POs G,H from 
Fig.\ref{al07h1}, and PO E from Fig.\ref{al3c2}, 
and {\em breathing and moving} vesicles, i.e., PO I from 
Fig.\ref{al07h1}, and POs F,G,H from Fig.\ref{al3c2}. 
Additionally, there are bifurcation from steady states 
to relative equilibria with non--zero translational Lagrange multipliers 
$\tau_x,\tau_y$, which corresponds to uniformly translating vesicles, 
see samples E,F in Fig.\ref{al07f1} and C in Fig.\ref{al3c1}. 
Naturally, HPs on such branches also yield breathing and moving 
vesicles, which moreover pick up some rotation. 

Next we used DNS to assess the stability of POs and more general 
dynamics. We distinguish between ``labframe'' DNS, or just DNS, 
and constrained DNS (cDNS), which expresses translations 
again via the multipliers $\tau_x$ and $\tau_y$.%
\footnote{We also ran cDNS with constraints for rotations and 
multiplier $\rho$, yielding the same results. Moreover, we ran cDNS 
in a simple implicit Euler version, and with the more elaborate 
high-order and adaptive code {\tt RADAU}, obtaining essentially 
identical results as long as solutions do not tend to self--intersection.} 
These two versions of DNS initially behave identical, but may 
deviate at large $t$: We believe that in (unconstrained) DNS numerical 
errors accumulate to excite the neutral (translational and/or rotational) 
modes, and as result for instance the   $D_4$ symmetric 
breather PO E from Fig.\ref{al3c2} 
is only metastable in DNS, but fully stable in cDNS, cf.~Fig.\ref{al3f}. 
Similarly, we believe that the $D_3$ symmetric (moving and rotating) 
breather F from Fig.\ref{al3c2} should be dynamically stable, 
but very quickly yields trigger waves under DNS, and also under 
cDNS is only metastable, leading to numerical failure after about 
20 periods.  On the other hand, DNS when not converging back to (stable) steady 
states generically goes to 
tumbling vesicles with periodic behavior of $\phi$ on $\ga$, 
i.e., source--sink pairs as in 
Fig.\ref{al3f}(a--d), Fig.\ref{al3f3}(a). 

We found analogous steady state bifurcations in 2D, and then focused on 
breathing vesicles, namely the two POs PO1 and PO2 from Fig.\ref{3f5a}.  
The general dynamics in 2D seem similar to 1D, but generally more 
complicated. 
In our experiments for $\al=2$ in \S\ref{3dsec}, when not starting 
near stable steady states, DNS at moderate $\zeta$ ($\zeta\approx 3.3$ 
in Fig.\ref{3f5b}(a)) yields 
tumbling vesicles with significant drifts, and typically between 4 and 
5 spikes (in $\phi$) on $X$. However, the tetrahedral breather PO1 
is again stable in cDNS, and cDNS at larger $\zeta$ ``typically'' 
yields convergence to (stable) steady states. 

As already said, we find the single morphogen $\phi$ 
mechanochemical model to be the simplest of its kind. 
Models that couple $X$ to several chemical species such as Brusselators, 
or reaction diffusion models involving  several different Min proteins, 
should be only slightly more complicated both theoretically and 
computationally, but naturally the parameter space becomes larger. 
Moreover, our motivating Fig.\ref{ef1} comes from experiments which 
call for bulk--surface models, in which oscillatory reactions in 
the fluid bulk play an important role. Such models can in principle also 
be studied in our setup, see, e.g., \cite[Ch.10]{p2pbook} for examples 
with fixed domains, but for dynamic domains 
the compatible motion of bulk and surface mesh points still needs 
further software development. 

In any case, we believe the bifurcation and continuation framework developed 
here (augmented by DNS for stability) to be crucial to get an overview of 
(relative) steady states and POs for the mechanochemical models \reff{mm1} 
and \reff{memo}, and similarly for the possibly more complicated models just 
indicated.

\appendix
\section{Some calculus}\label{lincal}
\subsection{Linearization}\label{lingen} 
For $\ga\subset\R^2$ a smooth, closed, parameterized (not
necessarily by arc-length) curve, and $F(\ga)$ an expression 
on some (Sobolev space) $H^r(\ga)$, we calculate the derivative via
\hueqst{ 
\pa_u F(\ga)u=\frac{\pa}{\pa\del} F(\ga+\del u
  \nu)\bigg|_{\del=0}, 
} 
with expressions 
evaluated at $\del=0$ if no argument is explicitly written. 
A key result is
the derivative of the metric $g=\spr{\pa_s\ga,\pa_s\ga}$ given by 
\hueqst{\pa_u g=2\spr{\pa_s (u\nu),\pa_s\ga}=2\sqrt{g}
  \spr{\pa_s\nu, \tau}u=-2g\kap u,
} 
using $\pa_s \nu=-\kap \sqrt{g}\tau$. 
 By the chain rule, we obtain
\hueqst{\pa_u\sqrt{g}=-\sqrt{g}\kap u.
}  
Since the Laplace--Beltrami
operator $\Del f=g^{-1/2}\pa_s(g^{-1/2}\pa_sf)$ depends on the
metric, its derivative is computed using 
$\pa_ug^{-1/2}=g^{-1/2}\kap u$ and the chain rule as 
\hualst{\pa_u\Del
  f&=
  \pa_u(g^{-1/2})\pa_sg^{-1/2}\pa_sf+g^{-1/2}\pa_u(\pa_sg^{-1/2}\pa_sf)\\
  &=\kap(\Del
  f) u + g^{-1/2}\pa_s (\kap g^{-1/2} u\pa_s f).
}  
Applying the product rule to the last term we arrive at 
\hueqst{\pa_u \Del
  f=2\kap (\Del f)u+\spr{\nabla (\kap u),\nabla f},
} 
where $\nabla f= (g^{-1/2} \pa_sf) \tau$. 
To linearize the curvature
$\kap=\spr{\Del \ga, \nu}$ we use the previous results and find
\hueqst{\pa_u \Del \ga= -(\Del u+\kap^2 u)\nu,
}
and consequently $\pa_u\kap=-(\Del u+\kap^2 u)$.

We can now write the full linearization of \reff{memo} in the normal
direction, as the remaining dependencies are simple polynomials. However, 
we first display our setup of the 4th order Helfrich flow in \reff{memo} 
as a 2--component 2nd order system.%
\footnote{Such a so called mixed formulation of a 4th order problems may be 
delicate, but seems 
unproblematic here because we have no boundary; see also \cite{SGJW22} for 
more general cases in the geometric setting.} With $\kap=\spr{\Del \ga,\nu}$ 
and $U=(u,\kap,\phi,\lam)$ 
we obtain 
\hual{\label{mmm}
\CM \pa_t U&=\bpm -\CG(U)\\ q_L(u)\epm, 
}
where $q(u)=L(\ga+u\nu)-L_0=0$, and 
\hugast{
 \CM=\bpm \frac{1+u}{g^{1/2}(u)}&&&\\ &0&&\\ &&1&\\&&&0
\epm,\quad 
-\CG(U)=\bpm -\Del(\kap{-}c_0)+\frac 1 2 \kap(\kap{+}c_0)(\kap{-}c_0)-\lam\kap\\
-\kap+\spr{\Delta\ga,\nu}u\\
D\Delta\phi-\del\phi+\zeta\frac{f(\kap)}{\om+f(\kap)}\epm. 
}
Thus, $\pa_{(u,\kap,\phi)} \CG\psi=  \pa_{u} \CG\psi_u+\pa_{\kappa} \CG\psi_\kappa+\pa_{\phi} \CG\psi_\phi$, 
where 
\hualst{ \pa_{u} \CG\psi_u=\bpm  -2\kap (\Del (\kap-c_0))\psi_u-\spr{\nabla (\kap \psi_u),\nabla (\kap-c_0)} \\
-\Del \psi_u-|v|^2 \psi_u \\ 
2\kap (\Del \phi)\psi_u+\spr{\nabla (\kap \psi_u),\nabla \phi}\epm,}

\hualst{ \pa_{\kappa} \CG\psi_\kappa=\bpm -\Del \psi_\kap-\frac 1 2(3\kap^2-c_0^2)\psi_\kap \\
 \psi_\kap \\
  \zeta \om \frac {\pa_\kap f}{(\om+f)^2}\psi_\kap \epm,\, \pa_{\phi} \CG\psi_\phi =\bpm  \beta(\Del \psi_\phi+\kap c_0\psi_\phi )\\
  0\\
    \Del \psi_\phi-\del\psi_\phi\\
  \epm,}
  and 
\hueq{
\label{fp} \pa_\kap f(\kap)=(1+\exp((\kap-\kap_0)\chi))^{-1}\exp((\kap-\kap_0)\chi).
} 
Importantly, the analytical Jacobian 
(instead of just numerical Jacobians) of $\CG$ 
is also used and very helpful in the numerics, 
for speed and accuracy. 

\subsection{Phase conditions}\label{pcs}
The equations \reff{nflow} and hence \reff{mmm} are purely intrinsic wrt the curve $\ga$, 
and hence normal perturbations that move $\ga$ in the
ambient space $\R^2$ are neutral directions, namely 
translations and rotations, aka rigid body motions. 
 For the numerical
continuation of steady states we have to take care of these
neutral direction, as they appear as eigenvectors to zero eigenvalues
of the Jacobian. 
Let $\Cm\in C^\infty(\R,\R^2)$ be a curve with $\Cm(0)=0$, and 
let $\vt\in C^\infty(\R,\R^{2\times 2})$ be a curve of rotation matrices 
with $\vt(0)=1$. They act on $\ga$ as 
\hueq{
T_{m(\Cr)}\ga{=}\{ \ga(s){+}\Cm(\Cr): s{\in}[0,2\pi)\} \quad\text{and}\quad 
R_{\vt(\Cr)}\ga{=}\{\vt(\Cr)\ga(s): s{\in}[0,2\pi)\}.
} 
$T_m$ and $R_\vt$ then form two dimensional respectively one dimensional 
Lie groups, as the tangent space at the
identity $\pa_\Cr|_{\Cr=0} T_{\Cm(\Cr)}\ga= \Cm'(0)$ is two dimensional, and
$\pa_\Cr|_{\Cr=0} R_{\vt}\ga= \bpm -\ga_2\\\ga_1\epm$ is one dimensional. 
As the Helfrich energy $\CH$ and the $\phi$ dynamics are 
intrinsically defined on $\ga$, the dynamics \reff{mmm} are 
equivariant under the group
generated by $\{T_m,R_\vt\}$. 
Letting 
\huga{ \label{pcon} 
q_1(u){=}\int_\ga \spr{e_1,u\, \nu}\dd s,\ 
q_2(u){=}\int_\ga \spr{e_2,u \,\nu}\dd s,\ 
q_3(u){=}\int_\ga \spr{\bpm -\ga_2\\\ga_1\epm, u\nu}\dd s
}
with the $\R^2$ scalar product $\spr{\cdot,\cdot}$, 
and introducing Lagrange multipliers $\tau_x$, $\tau_y$
and $\rho$ for the two positional and one rotational constraints, 
we add $\tau_x q_1(u)+\tau_yq_2(u)+\rho q_3(u)$ 
to  $\CH$. Because the constraints are linear in $u$,  we obtain 
the Euler Lagrange equation 
\hueq{\label{PCG} 
M\pa_t \bpm u\\\phi\epm 
=\CG(u,\phi)+\bpm \tau_x \spr{e_1,\nu}+\tau_y
  \spr{e_2,\nu}+\rho (\ga_1\nu_2-\ga_2\nu_1)\\ 0\epm
}
as the ``comoving frame equation'', 
together with the constraints $(q_1,q_2,q_3)(u)=0$, and the length 
constraint $q_L(u)=L(\ga_0+uN)-L_0=0$. Because all these constraints 
do not explicitly depend on the Lagrange multipliers, they all give 
DAEs of (differentiation) index 2. 
The time derivative of the constraints 
gives a linear algebraic system for the Lagrange multipliers 
\hueqst{\bpm 
\spr{\pa_uq_1,\pa_uq_1}_{L^2} & \spr{\pa_uq_1,\pa_uq_2}_{L^2} & \spr{\pa_uq_1,\pa_uq_3}_{L^2}\\
\spr{\pa_uq_2,\pa_uq_1}_{L^2} & \spr{\pa_uq_2,\pa_uq_2}_{L^2} & \spr{\pa_uq_2,\pa_uq_3}_{L^2}\\
\spr{\pa_uq_3,\pa_uq_1}_{L^2} & \spr{\pa_uq_3,\pa_uq_2}_{L^2} & \spr{\pa_uq_3,\pa_uq_3}_{L^2}\epm \bpm \tau_x\\\tau_y\\ \rho\epm= 
-\bpm \spr{\CG_1(U),\pa_uq_1}_{L^2}\\\spr{\CG_1(U),\pa_uq_2}_{L^2}\\
\spr{\CG_1(U),\pa_uq_3}_{L^2}\epm,
}
by using the chain rule and with for example 
$\spr{\pa_uq_1,\pa_uq_2}_{L^2}{=}\int_\ga \nu_1\nu_2\dd s $. Like for $\lam$ after 
\reff{nflow}, in principle this can be solved 
for $(\tau_x,\tau_y,\rho)$ 
as long $u$ is not too large, but here we just append $(q_1,q_2,q_3)$ to $q_L$ 
and extend the dynamical mass matrix $\CM$ from \reff{mmm} by three zero rows. 

In summary, for constant nonzero 
multipliers $\tau_x, \tau_y$ or $\rho$, a steady state of \reff{PCG} 
has some net rigid body motion in the lab-frame. Similarly, for 
POs $t\mapsto (\ga(t),\phi(t))$ we want the constraints \reff{pcon} 
(and $q_L(u)=0$) 
to hold at every $t$. Hence let
$\xi\in\Ga$ and 
$\ga=\xi_{-\eta(t)}\cdot \tilde \ga$ with a function $\eta(t)$, 
meaning a translation of $\tilde\ga$ in $e_1$
or $e_2$ direction, or a rotation by $-\eta(t)$. Then
\hueqst{ 
\spr{ \ga_t,\nu}= \spr{-\eta'(t) \pa_{\Cr}|_{\Cr=0}\xi_{\Cr}
    \cdot\tilde\ga+ \tilde\ga_t, \nu}=-\eta'(t) \spr{\pa_\Cr|_{\Cr=0}\xi_{\Cr}
    \cdot\tilde\ga, \nu} +\spr{\tilde\ga_t,\nu}.
}  
Dropping the
$\tilde\cdot$, the PCs for the dynamical problem can be considered
via \reff{PCG}, now with $\tau_x$, $\tau_y$ and $\rho$ time depending
as the derivatives of the corresponding $-\eta(t)$. 
Examples of nonzero periodic $\tau_x,\tau_y$ and $\rho$ are 
given in Fig.\ref{al07h1}(c), and Fig.\ref{al3c2}(b), 
which give relative POs, for which we then reconstruct the 
path of the center of mass in the lab-frame.  

All of these remarks hold analogously in 2D, where however we have 
to deal with three translational and two rotational rigid body motions.

\subsection{Amplitude equations}\label{aesec}
The linearization of \reff{locmemo} at $\MS^1$ with 
$\lam{=}\frac 1 2 (\al^2{-}1)$, i.e., $(u,\phi)=(0,0)$, 
yields $\pa_t U=\CL(\pa_s,\Lam)U$ with 
\hueq{
\CL(\pa_s,\Lam)= \bpm -(\pa_s+1)^2(\pa_s-1)^2 & \beta(\pa_s^2-\al)\\
		\frac{\zeta}{2\om} (\pa_s^2-1) & -D\pa_s^2-\del\epm, 
}
and the ansatz $U(s,t)=\exp(t\mu)e_m \Phi$ with $e_m=\exp(\ri ms)$, $m\in\Z$, 
and solving 
for $\re(\mu)=0$ yields \reff{zetc}, i.e., 
$\ds \zeta_c =\frac{2\om(m - 1)(m + 1)(m^2D^2 + \del)}{(m^2 - \al)\beta}$, 
with eigenvector 
\hueqst{\Psi_c= \bpm 1 \\ \frac {(m+1)^2(m-1)^2}{\beta(m^2-\al)}\epm e_{m}.
}
Here we briefly explain the amplitude formalism for branches 
bifurcating at $\zeta_c$, and how that depends on our smoothing 
parameter $\chi$ in \reff{smax}. 

Due to the $O(2)$ symmetry of the system, BPs on $\mS^1$ are double, 
but with the ``hidden $\Z_2$ symmetry'' of rotation by $\pi/(2m)$, which 
makes all bifurcations pitchforks. To derive the amplitude 
equations for these bifurcations, let 
$\zeta=\zeta_c+\eps^2\tilde \zeta$, introduce the slow time scale 
$T=\eps^2 t$, let  $A=A(T)$, and set 
\huga{\label{aea}
U(s,t)=\eps\Psi_A(s,t)= 
\eps A(T)\Psi_c(s)+\eps^2\left(\frac{1} 2\Phi_0e_0+\Phi_2 e_{2m}+\cc\right)+\hot. 
} 
We expand the right hand side of \reff{locmemo} as 
\hueq{ \CG(U)=\CL(\pa_s,\Lam) U+ B[U,U]+C[U,U,U]+\CO(\|U\|^4),}
with symmetric bilinear form 
$\ds B[U,U]=\frac 1 2 \left(\pa_u^2\CG+\pa_{u\phi}^2\CG+\pa_{\phi u}\CG+\pa_\phi^2\CG\right)[U,U]$, and symmetric trilinear form 
$C[U,U,U]=\frac 1 6\left(\pa^3_{u}\CG + 3\pa_{uu\phi}^3\CG+3\pa_{u\phi\phi}^3\CG +\pa_{\phi}^3\CG\right)[U,U,U]$, 
and similarly the constraint as 
\hueq{L(\ga)-L_0= \int_{\mS^1}\pa_u L(0) u\dd s + \int_{\mS^1} B_L(u,u)\dd s + \CO(\|u\|^4).}
We do not display here the formulas for $B$, $C$ and $B_L$ 
which become cumbersome, see \cite{mcSI}, but note that 
$(\pa_u^2\CG(u_1,u_2))_2=
 \frac{\zeta}{\omega ^2}(\omega \pa_\kap^2f(0)\pa^2_u\kap - 2\pa_uf(0)^2) u_2(s)u_1(s)$ 
 with $\pa_\kap f (\kap)$ from \reff{fp} and 
 \hual{\label{fpp} \pa_\kap^2f(\kap)=\frac{\chi\mathrm{e}^{\chi(\kap-\kap_0)}}{\left(\mathrm{e}^{\chi(\kap-\kap_0)} + 1\right)^{2}}.
}
Similarly, $\pa_\kap^k f(\kap)\sim \chi^{k-1}$, i.e., the higher 
derivatives become more and more singular with growing $\chi$. 

Inserting $\eps\Psi_A$ in \reff{locmemo}, all $O(\eps)$ terms vanish 
by construction, and at order $\eps^2$ we obtain 
\bsub\label{o2t}
\hual{
 0&=\CL(0,\Lam_c)\Phi_0+2|A|^2 \hat B[\Psi_c,\overline \Psi_c]+\tilde \lam \bpm 1 \\ 0\epm,\\
0&=\CL(i2m,\Lam_c)\Phi_2+A^2 \hat B[\Psi_c, \Psi_c],\\
0&= \pa_u L(0)(\Phi_0)_1+ |A|^2 \hat B_L[(\Psi_c)_1,(\overline \Psi_c)_1],  
}
\esub 
with $\CL(i2m,\Lam_c)$ from \reff{Lf}, and 
where $\hat B[\cdot,\cdot]= \CF B[\cdot,\cdot]$. From the dispersion relation for 
\reff{Lf} we can solve \reff{o2t} for $\Phi_0, \Phi_2\in\C^2$ and $\tilde\lam\in\R$. 

In order to remove the $\eps^3e_{m}$ terms we need to solve
\hualst{ \CL(im,\Lam_c)\tilde\Psi =& -\pa_TA\Psi_c + \pa_\zeta\mu_c\Psi_c+2A 
B[\Psi_c,\Phi_0]+2\overline A B[\overline \Psi_c,\Phi_2e_{2m}]\\
&+3 C[\Psi_c,\Psi_c,\Psi_c]
+\tilde \lam \pa_u \kap\bpm 1\\0\epm.}
As $\CL(m,\Lam_c)$ is not invertible, we get the solvability condition from the Fredholm alternative that 
\hueq{\label{ae1}
\pa_TA =\tilde \zeta \nu_1 A+\nu_2 A +\nu_3 \overline A +\nu_4 |A|^2A +A \tilde \lam \pa_u \kap ,}
where 
$\nu_1{=} \pa_\zeta\mu_c,\ \nu_2{=}2\spr{\hat B[\Psi_c,\Phi_0],\Phi_c^*}, \
\nu_3 {=}2\spr{\hat B[\overline \Psi_c,\Phi_2e_{2m}],\Phi_c^*}$, and \newline
$\nu_4{=} 3 \spr{\hat C[\Psi_c,\Psi_c,\overline{\Psi}_c],\Phi_c^*}$. 
Here  
$\hat C[\cdot, \cdot,\overline \cdot]= \CF C[\cdot, \cdot,\overline \cdot]$, and 
\hueqst{\Phi_c^*=\bpm 1 \\ \frac{(m^2 - \al)\beta}{(m^4 + (D - 2)m^2 + \del + 1)}\epm}
is the eigenvector to the adjoint eigenvalue problem $\CL^*(m,\Lam_c)\Phi_c^*=\overline \mu_c \Phi_c^*$, and so that $\spr{\Phi_c^*,\Phi_c}=1$. 
In \reff{ae1}, we explicitly kept $\nu_3$ because $\hat{B}$ involves 
$s$--derivatives and hence $\Phi_2e_{2m}$ must be used for $\nu_3$, not just 
$\Phi_2$. However, as 
$\nu_2, \tilde \lam\sim |A|^2$ and $\nu_3\sim A^2$ we get the normal 
form amplitude equation for a pitchfork 
\hual{\label{fcgl} \pa_T A= \tilde \zeta\nu_1 A +\tilde \nu |A|^2 A.} 

\begin{figure}[ht]
\btab{l}{{\sm (a1)} \hs{40mm}{\sm (a2)}\\
\ig[width=0.27\tew]{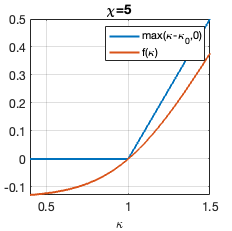}\hs{5mm}\ig[width=0.35\tew]{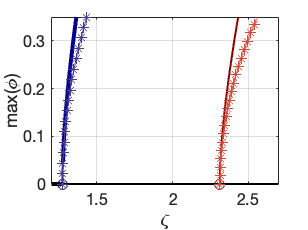} \hs{0mm}\ig[width=0.35\tew]{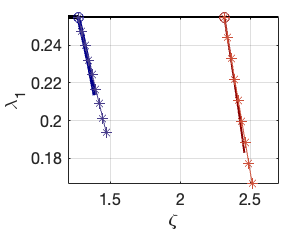}\\
{\sm (b1)} \hs{40mm}{\sm (b2)}\\
\ig[width=0.27\tew]{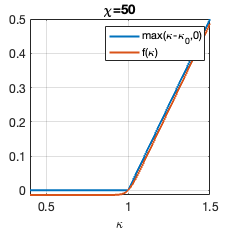}\hs{5mm}\ig[width=0.35\tew]{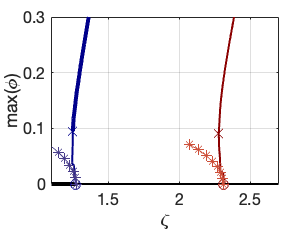} \hs{0mm}\ig[width=0.35\tew]{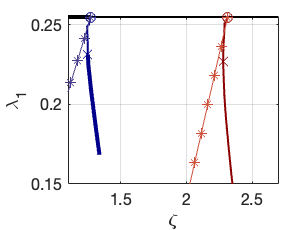}
}

\vs{-2mm}
\caption{{\sm  Comparison of the numerical continuation (lines) and the amplitude prediction (dotted lines) for $\chi=5$ in (a) and $\chi=50$ in (b) for the parameters $(\al,\beta,D,\del, \om,)=(0.7,1,0.1,0.3,1)$.
\label{ampf1}}}
\end{figure}

To find steady nontrivial steady state branches bifurcating at $\zeta_c$, 
we can now solve \reff{fcgl} for steady states $A$. 
However, the quality of the prediction by the amplitude equation
heavily depends on the value of the smoothing parameter $\chi$. 
From \reff{fpp}, the
quadratic terms scale as $\chi$ and the cubic 
terms as $\chi^2$. This couples via $\nu_2$
and $\nu_3$ in the cubic coefficient in the amplitude equation. For
the parameters $(\al,\beta,D,\del, \om)=(0.7,1,0.1,0.3,1)$ we compare
in Fig.\ref{ampf1} the predictions of the amplitude equations with
the numerical continuation for different values of $\chi$. For small
$\chi$, the curvature of $f_\chi$ 
at $\kap=1$ is small and the approximation of the
$\max$ function rather poor. For $\chi=5$, 
the numerical values in \reff{fcgl} are 
\bsub\label{ae2}
\hual{
\CW_2:\quad \pa_TA&=0.51 A\tilde \zeta - 19.85|A|^2A,\\
  \CW_3:\quad \pa_TA&=0.51 A\tilde \zeta - 244.33|A|^2A. 
} 
\esub 
In this case, 
the bifurcations are supercritical and the predictions by the
amplitude equations for $\phi$ and also for the Lagrange multiplier are
rather accurate also at $\CO(1)$ amplitude. In contrast, for higher
$\chi=50$, the amplitude equations are 
\bsub\label{ae3}
\hual{
\CW_2:\quad \pa_TA&=0.51 A\tilde \zeta + 583.4|A|^2A,\\
  \CW_3:\quad \pa_TA&=0.51 A\tilde \zeta + 5708.38|A|^2A. 
} 
\esub 
From \reff{ae3}, the bifurcations are now subcritical. 
This also holds for our numerical bifurcations, 
but with a fold shortly after the
bifurcation point, and for large $\chi$  the approximation by the AEs 
is therefore only good close to bifurcation. 

In summary, for any $\chi>0$, the AEs \reff{fcgl} can be justified by 
standard center manifold theory, but the range of validity, i.e., 
the maximal allowed $\eps$ in \reff{aea} shrinks with $\chi\to\infty$.  
However,  
the choice of $\chi$ does not seem to have much influence 
on the global behavior of the nontrivial steady state branches $\CW_m$; 
this is expected as for nontrivial solutions, $\kap$ is bounded 
away from 1 on most of $\ga$. 
A rather ambitious but interesting next step would be to study the
local bifurcations from $\mS^1$ in the singular limit $\chi\to\infty$, 
where standard bifurcation theory does not apply. See also Remark \ref{smrem}. 

\section{Numerical algorithms}\label{algo}
The {\tt Xcont} extension \cite{geomtut,geompap} of \pdep\ provides 
methods for the numerical continuation and bifurcation analysis 
of manifolds as solutions of geometric PDEs. The focus so far 
were steady 2D problems, but the algorithms also apply to the 1D case, 
see the introduction of \cite{AMPHD}. Some further features needed 
here pertain to PO continuation, explained below, after first briefly 
recalling the geometry setting of the {\tt Xcont} extension, in 1D. 
\subsection{Spatial and temporal discretization}\label{timevo}
We spatially discretize a closed curve $\ga$ into ${\tt n_t}$ linear
connections $({\tt tri})_{i=1}^{\tt n_t}$ with ${\tt n_p}={\tt n_t}$ nodes
$({\tt g}_j)_{j=1}^{\tt n_p}$, ordered such that
the two neighbors of ${\tt g_j}$ are ${\tt g_{i-1}}$ and
${\tt g_{i+1}}$ (with ${\tt g_0}={\tt g_{np}}$ and ${\tt g_{np+1}}=
{\tt g_1}$ implied). Using piecewise linear hat functions $\psi_j$, with
$\psi_j({\tt g}_k) = \delta_{jk}$, the weak form of the
Laplace–Beltrami operator is given by 
\hueqst{ L_{jk} = -\int_\ga
  \spr{\nabla \psi_j, \nabla \psi_k}\dd s, 
} 
where $\nabla f = (\pa_s f )\tau$ is the curve gradient. Instead 
of the standard mass matrix $M_{jk} = \int_X \psi_j \psi_k \dd S$, 
as in \cite{geompap} we use the (diagonal, lumped) Voronoi 
mass matrix, which in 1D simply reads 
\hueq{\label{vmass} M_{jj} = \frac
  {1}{2} \left(\| {\tt g_{j}}-{\tt g_{j-1}}\|_2+\|{\tt g_{j+1}}-{\tt
    g_{j}}\|_2\right).
} 
Then we discretize $\kap$ in the weak FEM
formulation as 
$M {\tt k} =\spr{L {\tt g},\nu}$, where $\nu$ is 
the vertex normal computed as the average of
the two adjacent face normals,  weighted by element size. 
The implementation uses the {\tt gptoolbox} \cite{gpgit}, 
see also \cite{geompap} for comments
on convergence rates (in 2D). Finally, we 
discretize the 
mixed formulation \reff{mmm} 
with a singular mass matrix, and discretize the reaction diffusion 
equation of the morphogen 
in a standard fashion, i.e., 
\bsub\label{sdflow}  
\hueq{
\bpm M & \\ & 0\\ & & M\epm \bpm \spr{ \dot {\tt g},\nu}\\ \dot {\tt k} \\
  \dot \phi \epm 
{=} - \tilde {\tt G}(u,{\tt k},\phi)
{=} -\bpm L ({\tt k}-c_0){+}M(\frac 1 2 {\tt k}\big({\tt k}
{+}c_0\big)\big({\tt k}{-}c_0\big) {-}\lam{\tt k})\\
  \spr{L {\tt g},\nu}{-}M{\tt k} \\
  -L \phi +M\left(\del \phi -\zeta \frac{ f({\tt k})}{\om +f({\tt
      k})}\right)\epm,}  
together with the (discretized) length constraint 
\huga{
q_L(u)=L(\tt g)-L_0=0,} 
\esub
which determines the Lagrange multiplier $\lam$, and similar discretizations 
of the phase conditions $(q_1,q_2,q_3)(u)=0$, and with then adding (the 
discretizations of) 
$\tau_x\spr{e_1,\nu}+\tau_y
  \spr{e_2,\nu}+\rho (\ga_1\nu_2-\ga_2\nu_1)$ from \reff{PCG} to the first line 
of the rhs of (\ref{sdflow}a). 

\paragraph{Steady state continuation.}
Steady states can now be continued in parameters 
as solutions of 
\huga{\label{sdflowa}
0=-\tilde {\tt G}(u,{\tt k},\phi;\Lam),\qquad q_L(u)=0, 
}
where $\Lam\in\R^p$ is a generic name for the pertinent {\em active} parameters, 
e.g., $\Lam=(\zeta,\lam)$ with $\zeta$ the {\em primary active} parameter here, 
and $\lam$ as secondary active parameter. However, due 
to the translational and rotational invariance of \reff{memo} and hence 
approximate%
\footnote{the discretization already breaks these invariances, but only 
very weakly and in an uncontrolled way} 
 invariances of \reff{sdflow0}, as explained above we need phase conditions 
$\qp(\ga,\phi,\Lam)=0$ 
and associated Lagrange multipliers as further 
secondary active parameters, see \S\ref{pcs}, and for the system 
thus obtained we write 
\huga{\label{sdflow0}
0=-\tilde {\tt G}(u,{\tt k},\phi;\Lam),\qquad q(u)=0. 
}
Now assume we have a point 
\huga{\label{uuu}
U_0=(\ga(\sig_0),\kap(\sig_0),\phi(\sig_0),\Lam(\sig_0))
}
on a solution branch $(\sig_-,\sig_0]\ni\sig\mapsto U(\sig)$ 
(with a slight abuse of notation explained below 
writing $\ga(\sig_0)$ instead of $u(\sig_0)$ for the first component of $U$)   
and a unit (in some suitable weighted norm) 
tangent 
\huga{
\Tau(\sig_0)=\pa_\sig U(\sig_0)=(\Tau_\ga,\Tau_\kap,\Tau_\phi,\Tau_\Lam)
}
to that solution branch, 
and a step length $\ddsi$.%
\footnote{The default name for the (dummy) arclength parameter 
in continuation problems is $s$, see \cite{p2pbook}, also in our 
{\tt Xcont} setting \cite{geompap}; however, for the 1D problem we 
already use $s$ as the parametrization of $\ga$. Also note that 
the branch $\sig\mapsto U(\sig)$ always means a branch of curves $\ga$, 
fields $\phi$ and (active) parameters $\Lam$.} 
Then we make a predictor 
\huga{
U(\sig)=U(\sig_0)+\ddsi \Tau(\sig_0) 
} 
at $\sig=\sig_0+\ddsi$, and aim to solve \reff{sdflow0} for 
$(u,\kap,\phi,\Lam)$ by Newton loops, starting with $u^0=0$, 
in the hyperplane orthogonal to $\Tau$, which is the 
crucial idea of arclength continuation to deal with folds. 
In case of success we update $\ga(\sig)=\ga(\sig_0)+u(\sig)\nu(\sig_0)$ and 
go to the next step, which we summarize in short as {\tt updX} 
(as the name of the pertinent function in {\tt Xcont}, ``update $X$'', in 1D  
$X=\ga$); otherwise, a standard idea is to reduce the step length 
$\ddsi$ and try again, and altogether the convergence speed (or failure) 
of the Newton loops is a basis for stepsize control. Thus, $u$ is always 
small in this setup and only meaningful together with the ``current'' base 
manifold $\ga(\sig_0)$, which is why we write $\ga(\sig_0)$ in \reff{uuu}.

\paragraph{PO continuation.}
To continue POs, we also need a temporal discretization
of \reff{sdflow}.  To compute a PO with (unknown) period $T$, 
we consider the boundary value problem rescaled to the temporal
interval $[0,1]$, i.e., 
\bsub\label{pon}
\hueq{
 \bpm M & \\ & 0\\ & & M\epm \bpm \spr{ \dot {\tt g},\nu}\\ \dot {\tt k} \\
  \dot \phi \epm + T \tilde {\tt G}(u,{\tt k},\phi)=0, 
}
together with the periodicity condition 
\huga{
(u(0),{\tt k}(0),\phi(0),\Lam(0))=(u(1),{\tt k}(1), \phi(1),\Lam(1)), 
}
\esub 
where we already preview that at least some of the Lagrange-multipliers,  
namely $\lam$, are dynamic, i.e., functions of $t$. 
For  a discretization $0=t_0,t_1,\ldots,t_{m-1}\in[0,1)$ with stepsize $h_j=t_{j+1}-t_j$ and periodicity condition $U(t_m)=U(t_0)$, 
we discretize the temporal derivative of the Helfrich
equation as
$\spr{\dot{\tt g}(t_j),\nu(t_j)}\approx \spr{h_j^{-1}({\tt g}(t_j)-{\tt
    g}(t_{j-1})),\nu(t_j)}$, and use a mid
point approximation of the first and third component of $\tilde{G}$, i.e., 
$\tilde G_i(t_j)=\frac 1 2(\tilde G(t_j)+\tilde G(t_{j-1}))$, but 
evaluate the second component for the elliptic equation $\spr{Lg,\nu}-Mk=0$, 
and the constraints, at $t_j$. 
In the numerical continuation of POs, the updating procedure of $\ga$ 
is as follows: The current $\ga$ is given 
by the field ${\tt p.hopf.X}$, and the new $\ga$ at time
slice $t_j$ is 
${\tt Xn}(j)= {\tt p.hopf.X}(:,:,j)+ {\tt u}(:,j)\nu(j)$. The
temporal discretization of $\spr{ \dot {\tt g},\nu}$ in the first line of (\ref{pon}a) then
 is 
\huga{\label{ddd1}
\spr{\dot{\tt g}(t_j),\nu(t_j)}\approx \spr{h_j^{-1}({\tt Xn}(t_j)-{\tt Xn}(t_{j-1})),\nu(t_j)}.
}

\paragraph{DNS and cDNS.} 
In DNS, we do one time step of length $h$ 
and then update $\ga$, and hence \reff{ddd1} simplifies to
\hual{ \spr{\dot{\tt g}(t+h),\nu(t+h)}&\approx 
\spr{h^{-1}({\tt g}(t+h)-{\tt g}(t)),\nu(t+h)}\notag\\
  &= \spr{h^{-1}({\tt g}(t)+u\nu(t)-{\tt g}(t)),\nu(t+h)}=h^{-1}\spr{\nu(t),\nu(t+h)}u.\label{ddd2}
}
Due to the length constraint $q_L(u)=0$ we have a differential 
algebraic equation (DAE) of (differentiation) 
index 2, where the Lagrange multiplier $\lam$ is a 
time dependent part of the solution. With $\eta(t):=\spr{\nu(t{-}h),\nu(t)}$ 
we write the full DAE system as 
\hueq{\label{dae1}
 \bpm M\eta(t) & \\ & 0\\ & & M\\ & & & 0\epm \bpm \dot u(t) \\ 
\dot {\tt k}(t) \\ \dot \phi(t) \\ \dot\lam(t) \epm + 
\bpm \tilde {\tt G}_1(u(t),{\tt k}(t),\phi(t),\lam(t))\\ 
\tilde {\tt G}_2(u(t),{\tt k}(t),\phi(t),\lam(t))\\ 
\tilde {\tt G}_3(u(t),{\tt k}(t),\phi(t),\lam(t))\\ 
q_L(u(t))\epm=\bpm 0\\ 0\\ 0\\0 \epm, 
}
where again the $(\dot\,)$ stands for the discretization of the time 
derivative, and $\dot\lam$ and $\dot{\tt k}$ only appear formally 
in \reff{dae1} to visualize the structure. 
We call \reff{dae1} with $\tau_x,\tau_y,\rho=0$ 
the lab--frame formulation of the DAE 
(where the 2nd row of \reff{dae1} has differentiation index 1 and the last row index 2), 
and the associated numerical integration simply DNS, see below. 

Alternatively, including the positional constraints 
$q_1,q_2,q_3$ from  \reff{pcon},\reff{PCG}  in the DNS 
yields DAEs for $(\tau_x,\tau_y,\rho)$ of index 2. We write this 
as 
\hueq{\label{dae2}
 \bpm M\eta(t)  \\ & 0\\ & & M\\ & & & 0\\ & & & & 0_p\\
\epm \bpm \dot u(t) \\ 
\dot {\tt k}(t) \\
  \dot \phi(t) \\ \dot \lam(t)\\\dot\Lam_p \epm + 
\bpm \tilde {\tt G}_1(u(t),{\tt k}(t),\phi(t),\Lam(t))\\
 \tilde {\tt G}_2(u(t),{\tt k}(t),\phi(t),\Lam(t))\\
\tilde {\tt G}_3(u(t),{\tt k}(t),\phi(t),\Lam(t))\\
q_L(u(t))\\ q_p(u(t))\epm=\bpm 0\\0\\0\\0\\0\epm 
}
where $q_p(u(t))=(q_1,q_2,q_3)(u(t))$, $0_p=0\in\R^{3\times 3}$, 
and $\Lam_p=(\tau_x,\tau_y,\rho)$, 
or, dropping the rotational constraint $q_p(u(t))=(q_1,q_2)(u(t))$,  
$0_p=0\in\R^{2\times 2}$, 
and $\Lam_p=(\tau_x,\tau_y)$, and we call \reff{dae2} 
constrained DNS (cDNS). 
The idea of using \reff{dae2} instead of \reff{dae1} for DNS is 
to trade drift (and rotation) in the lab frame for time--dependent non--zero 
$\tau_x,\tau_y$ (and $\rho$),  which turns out to be a more robust 
formulation as it seems to avoid the accumulation of numerical errors 
in the rigid motion neutral modes. 

There are several numerical strategies for DAEs, see for example 
\cite{hlr89,hw2}. A simple method for \reff{dae1} and \reff{dae2} 
is implicit Euler 
such that the 2nd and 4th row in \reff{dae1} (and the 5th row in 
\reff{dae2}) simply yield constraints at $t_{j+1}$, and we 
take this as a our basic strategy. However, 
 implicit Euler is only first order in $h=t_{j+1}-t_j$, and 
can be ``overdamping'' and hence in particular be inadequate for 
approximating POs. See, e.g., the discussion 
of the mathematical pendulum in Cartesian coordinates as (index 1,2, or 3) 
DAE in \cite{hlr89}. Therefore we checked our implicit Euler results 
for \reff{dae1} and \reff{dae2} against the high--order adaptive code 
{\tt RADAU}. The results fully agree 
for stable POs and for the convergence of DNS to stable steady states, 
but for unstable dynamics (trigger waves, and cases that lead 
to self--intersection and blowup) significant differences can 
occur, as {\tt RADAU} often yields non--convergence 
as $\ga$ becomes strongly distorted. Thus, we mostly use {\tt RADAU} 
for validation of the stable dynamics. Moreover, 
the 2D results in \S\ref{3dsec} are all with implicit Euler as 
{\tt RADAU} becomes too slow for the many degrees of freedom in 2D, 
and is only used for some validations. 

\subsection{Mesh handling}\label{meshh}
\begin{figure}[ht]
\btab{lll}{
{\small (a)} \hs{45mm}{\small (b)}\vs{2mm}\\
%
\hs{0mm} \ig[width=0.35\tew]{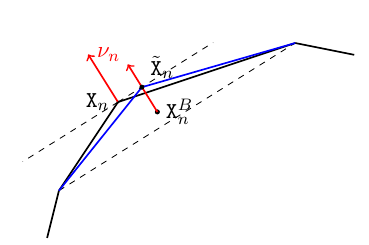} 
\hs{5mm} \ig[width=0.3\tew]{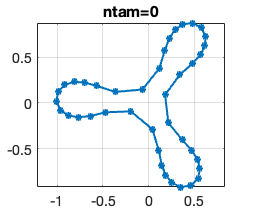} \hs{2mm}\ig[width=0.3\tew]{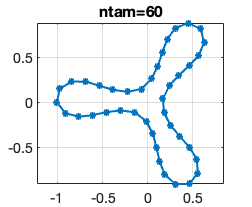}
}
\caption{{\small (a) Schematic picture of the {\tt tamo} algorithm for constructing a the new vertex $\tilde {\tt X}_n={\tt X}_n^B+\eta\nu_n$. (b) shows an exemplary use of ntam times application of {\tt tamo.m}}\label{tamof1}}
\end{figure}

Mesh handling is crucial in all our numerics, i.e., steady state continuation, 
PO continuation, 
and DNS: as $\ga$ changes (in time, or in continuation), 
the mesh may deteriorate, in particular since so far we only describe normal 
displacements. Basic refinement and coarsening strategies for
geometric PDE's in \pdep\ are explained in \cite{geomtut} 
(for the 2D case). Especially
coarsening is crucial to deal with neck
development. Here we additionally rely on a different type of mesh
handling following \cite{BMN05}: {\tt tamo}=tangential mesh optimization, 
see Fig.\ref{tamof1}(a) for illustration, and (b) for examples on a 
coarse mesh.

The idea is to find a good 
approximation of the tangent space to move a point along, 
keeping the enclosed area fixed. We use 
the tangent plane (in 1D tangent line) spanned by adjacent barycenters. 
Let ${\tt X}_n$ be a vertex of the polygonal curve and $\CS$ the ``star''
surrounding it (the support of the FEM basis
function). In 1D, $\CS$ is simply the two elements connecting
${\tt X}_{n-1}$ with ${\tt X}_{n}$, and ${\tt X}_{n}$ with
${\tt X}_{n+1}$. Compute the vertex normal $\nu_n$ as the weighted
average of each element normal, 
\hueqst{\nu_n
  =\frac1 {|S_{n-1}|+|S_n|}
  (\nu_{S_{n-1}}|S_{n-1}|+\nu_{S_{n}}|S_{n}|). 
}  
The barycenter of an
element is $\hat {\tt X}_n{=}\frac 1 2 ( {\tt X}_n+ {\tt X}_{n-1})$ and
the averaged barycenter of a star then is 
$\hat{\tt
    X}^B_n{=}\frac 1 2 ( \hat {\tt X}_n +\hat {\tt X}_{n-1})$.  
From this we compute the new position $\tilde {\tt X}_n$ as
\hueqst{\tilde {\tt X}_n=\hat{\tt X}^B_n+\eta \nu_n,
} 
where $\eta$ is chosen such that the areas of triangles 
$({\tt X}_n,\hat {\tt X}_n,\hat {\tt X}_{n-1})$ and
$(\tilde{\tt X}_n,\hat {\tt X}_n,\hat {\tt X}_{n-1})$ are equal. 
Using that the area of a triangle 
$(A,B,C)$ is 
$\frac 1 2 (B-A)\times (C-A)$,  we thus want
$(\tilde {\tt X}_n- \hat {\tt X}_n)\times (\hat{\tt
    X}_{n-1}-\hat{\tt X}_n)=({\tt X}_n- \hat {\tt X}_n)\times
  (\hat{\tt X}_{n-1}-\hat{\tt X}_n)$,  
and solving for $\eta$ yields 
\hueqst{\eta = \frac{({\tt X}_n- \hat {\tt X}_n)\times (\hat{\tt
      X}_{n-1}-\hat{\tt X}_n) -(\hat {\tt X}^B_n- \hat {\tt
      X}_n)\times (\hat{\tt X}_{n-1}-\hat{\tt X}_n)}{ (\nu_n- \hat
    {\tt X}_n)\times (\hat{\tt X}_{n-1}-\hat{\tt X}_n)}.
} 
This algorithm {\tt tamo} 
is fast and under some iteration ({\tt ntam} iterations with 
typically {\tt ntam} between 10 and 80) very efficiently improves meshes. 
Essentially the same works in 2D, in a slightly more complex fashion, 
and we apply {\tt tamo} as follows: 
\bci 
\item In steady state continuation, we apply {\tt tamo} after each 
continuation step. Importantly, as $\phi$ lives at the mesh points, 
after each {\tt tamo} step, or alternatively at the end of 
the {\tt tamo} iteration, $\phi$ must be interpolated from 
the old mesh points to the new mesh points, and we call 
this {\tt phii}. This is very 
easy in 1D, by simply using the arclength coordinates on $\ga$
 and $\ga_{\text{new}}$, but slightly more complicated and expensive in 2D. 
Then, after {\tt tamo} and {\tt phii}  we run one more Newton loop to update 
the current solution $(\ga,\phi)$, and then proceed to the 
next continuation step. 
\item In PO continuation, for larger amplitude POs with 
stronger distortion of $\ga$, 
we apply {\tt tamo} (and {\tt phii}) at each time 
slice $t_j$, and then run a Newton loop on the full system \reff{sdflow}. 
\item In DNS, we run {\tt tamo} and {\tt phii} after a selected 
number of time steps ${\tt n_{DNS}}$ of step size $h$. 
For slow dynamics, e.g., near 
steady state, ${\tt n_{DNS}}$ (and $h$) can in principle be chosen large, 
but in faster dynamics 
choosing small ${\tt n_{DNS}}$, in some cases ${\tt n_{DNS}}=1$, 
is very helpful to have DNS with $h$ bounded away from 0. 
\eci 

\subsection{Numerical details, and mesh refinement}
\label{ncsec}
\def\res{\|G(u)\|_\infty}\def\tol{{\tt tol}}
We close with listing the numerical settings used 
in the main text, some experiments with mesh refinement in space and 
time, and, following Fig.\ref{3f6} some examples of 2D vesicles at 
larger $\zeta$. 

The numerics for continuation use $\res=\tol=10^{-8}$ throughout, where 
for steady states $G(U)$ stands for \reff{sdflow0}, 
while for PO continuation  $G(U)$ stands for \reff{pon}. 
In 1D, we initially discretize $\mS^1$ into $300$ equal elements, yielding 
$900$ DoF, while for the $t$--discretization of POs 
we use $60$ to $100$ points, yielding up to $90\,000$ DoF. 
The initial discretization of the sphere is by subdivision and 
projection, leading to $1280$ triangular elements, and 
$1926$ spatial DoF, for POs combined with 30 points in time, 
hence approximately $60\,000$ total DoF.%
\footnote{However, even with a comparable number of DoF in 1D and 2D, 
the 1D PO computations are still much faster, because they can 
better exploit bordered linear system solvers: in 1D, the spatial 
discretization essentially yields 
tridiagonal matrices at each time slice, which can efficiently 
be dealt with by LU decomposition, which becomes much more expensive in 2D.} 
In summary, with the above discretizations, 
computation time for 1D POs is on the order of 30s per step 
on a standard desktop PC, while for 2D POs it can easily be on the order 
of 10 minutes per step. For steady state continuation, computation time 
is not a problem ($<1$s per step in 1D, $<10$s in 2D), and also 
the implicit Euler method for DNS is fast (again 
$<1$s per step in 1D, $<5$s in 2D), while {\tt RADAU} is significantly 
slower. Thus, the main bottleneck, and main reason for our rather 
sparse discretizations in 2D, are the POs in 2D. 
For both, 1D and 2D, we generally checked that finer spatial 
or/and temporal resolutions yield no qualitative changes; 
see Fig.13 for one example of mesh--refinement on the sphere 
yielding a reduction of the (mesh--asymmetry--induced) 
translational multipliers from $10^{-4}$ to $10^{-5}$ (and smaller), 
and below we report some further such experiments on mesh refinement. 

For the DNS we generally use $\tol{=}10^{-12}$ for our 
implicit Newton method, and adaptive time--stepping (based on 
the convergence speeds of the Newton loops) with stepsize $h$ on 
the order of $h{=}10^{-4}$ to $h{=}0.1$. For {\tt RADAU} we 
use the standard settings with ${\tt atol}{=}10^{-5}$, 
${\tt rtol}{=}10^{-7}$, which leads to stepsizes between 
$h_{\text{min}}{=}10^{-7}$ and $h_{\text{max}}=0.2$. 

\begin{figure}[ht]
\btab{ll}{{\sm (a)}&{\sm (b)}\\[-0mm]
\rb{0mm}{\ig[width=0.26\tew]{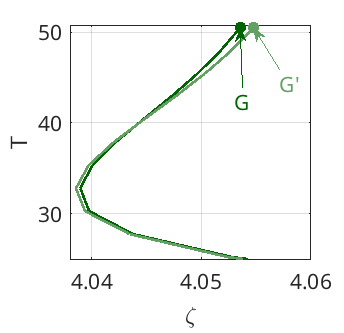}\hs{-1mm}
\ig[width=0.25\tew]{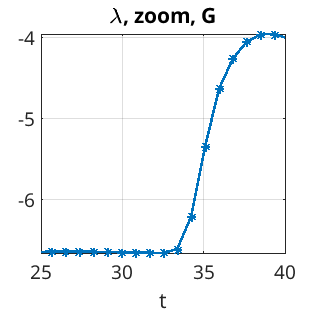}\hs{-1mm}
\ig[width=0.25\tew]{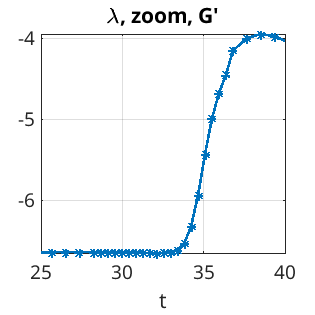}}
&\hs{-3mm}\ig[width=0.26\tew]{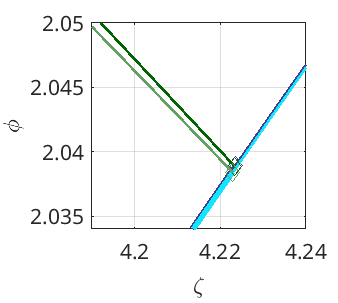}
}
\vs{-2mm}
\caption{{\sm  (a) 
Comparison of the PO branch (dark green) with sample G from 
Fig.\ref{al07f1} and Fig.\ref{al07h1} with the same 
branch (light green) obtained after temporal mesh--refinement 
near the kink in $\lam$. (b) Comparison of steady wave number 2 
branch from Fig.\ref{al07f1} 
(dark blue) with the same branch after mesh--refinement from $300$ to $360$ 
elements (light blue), and the bifurcating PO branches (dark/light) green. 
\label{nuf1}}}
\end{figure}

In Fig.\ref{nuf1}  we report experiments with temporal and spatial mesh 
refinement following Fig.\,\ref{al07f1}, specifically pertaining to HP1 
on the spatial wavenumber $2$ branch, and 
the PO bifurcating from this with sample G in Fig.\ref{al07f1}, 
for two reasons: Adaptive spatial mesh 
refinement is more significant on solutions more strongly distorted 
from the circle, and the PO branch with G shows relaxation 
oscillations which call for mesh adaptation in $t$. 
In (a), after 5 continuation steps (near $\zeta=4.131$) 
we use adaptive temporal mesh refinement near 
the kink visible in $\lam$ to (rather mildly) refine from 
$m=60$ to $m=75$, yielding the light green branch with sample G' 
with comparison to the original dark green branch with sample G. 
While there are some quantitative differences in the displayed period $T$, 
and similarly in other diagnostics such as $\spr{\phi}$, no 
qualitative changes occur. 
In (b) we do a purely spatial mesh--refinement on the steady state branch 
at the two tips of high curvature, going from $n_t=300$ to $n_t=360$ elements. 
The difference on the steady state branch then is rather small, e.g., 
$\zeta_{\text{HP}}=4.2234$ for $n_t=300$ and $\zeta_{\text{HP}}=4.2230$ for $n_t=360$. 
Near bifurcation there then also is a notable difference in $\spr{\phi}$ 
for the bifurcating PO branches, 
dark green with $n_e=300$ and light green with $n_e=360$, both 
with $m=60$ temporal discretization points, but their 
overall behavior is again very similar, like in (a). 

Similar conclusions hold for all other 1D results which we tested via 
(spatial or temporal, or both) mesh refinement, and in summary we find 
that our admittedly somewhat coarse meshes correctly resolve the 
qualitative features. Moreover, we like to stress again that for 
strongly deformed $\ga$ (or $X$ in 2D), the tangential mesh--optimization 
{\tt tamo} (see Fig.\ref{tamof1}(b) for a cartoon picture with 
{\em very} coarse 
resolution) seems more important and efficient for robust results, and hence 
is used throughout. 

\begin{figure}[ht]
\btab{l}{{\sm (a)}\\[-0mm]
\hs{-3mm}\ig[width=0.5\tew]{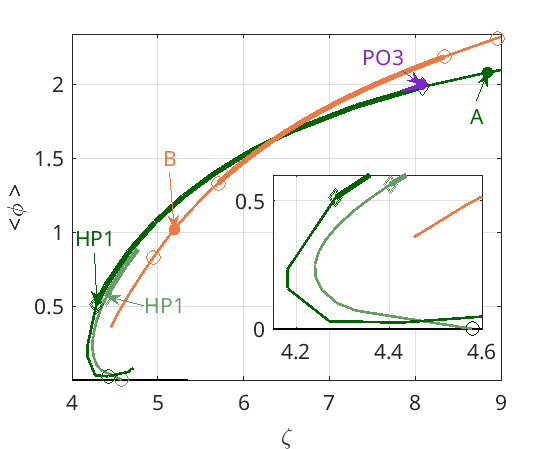}
\rb{28mm}{\btab{l}{\ig[width=0.23\tew]{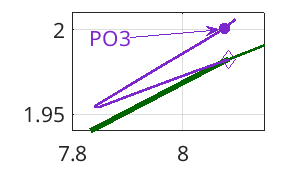}\\
\ig[width=0.27\tew]{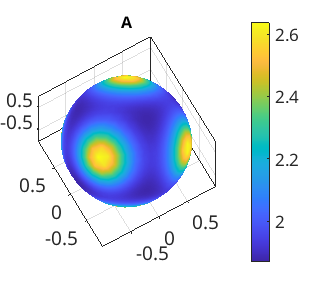}\hs{-1mm}
\ig[width=0.27\tew]{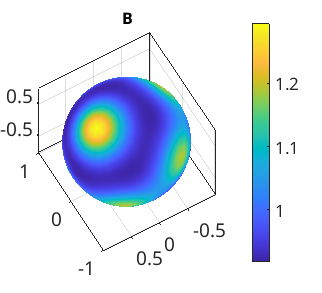}}}\\
{\sm (b)}\\[-3mm]
\ig[width=0.23\tew]{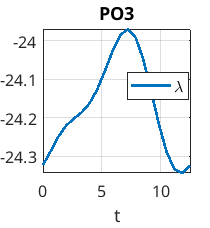}\hs{3mm}
\hs{-2mm}\ig[width=0.22\tew,height=40mm]{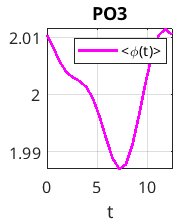}\hs{3mm}
\ig[width=0.27\tew]{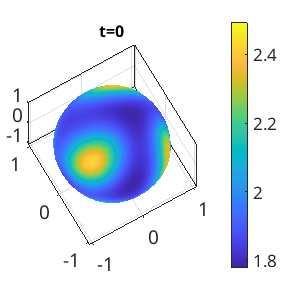}\hs{-1mm}
\ig[width=0.27\tew]{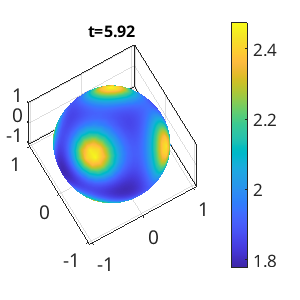}}
\vs{-2mm}
\caption{{\sm  (a) Continuation of octahedral (dark green) and decahedral 
(orange) branches on refined meshes from Fig.\ref{3f6}, with respective 
starting point from the DNS at $\zeta{=}6$ and $\zeta{=}7$. 
Light green is the octahedral branch from Fig.\ref{3f5a}. 
(b) sample from PO3 branch at $\zeta{\approx}8.075$. 
\label{nuf2}}}
\end{figure}

In Fig.\ref{nuf2} we present an octahedral  
(dark green, $n_t{=}2848$ elements) and a decahedral branch 
(orange, $n_t{=}2992$ elements) obtained from starting 
at the end of the cDNS in Fig.\ref{3f6}(a) at $\zeta=6$ and 
$\zeta=7$, on the refined meshes. We compare the octahedral branch 
to the 
``original'' octahedral branch from Fig.\ref{3f5a} (light 
green, $n_t=1280$ elements). 
Going to smaller $\zeta$ on the dark green branch 
we find the same qualitative behavior 
as on the light green branch, with a HP at $\zeta\approx 4.285$ 
($\zeta\approx 4.402$ on light green), and 
the fold at $\zeta_{\text{fold}}\approx 4.182$ 
($\zeta_{\text{fold}}\approx 4.241$ on light green). 
The zoom inset also shows that the dark green branch then 
directly 
continues past the shere to another branch from the BP of multiplicity 9 
on the sphere at $\zeta\approx 4.579$, but this is the natural 
behavior in the backward continuation. 
Thus, while the locations of the FPs and the HPs at low $\zeta$ 
differ perceivably  between dark and light green, 
the branches are qualitatively the same, 
and a PO analogous to  PO2 from Fig.\ref{3f5a} also 
bifurcates from the refined branch at low $\zeta$ (not shown). 

Next, we also follow the octahedral branch 
to large $\zeta$ and find that it loses stability in another HP 
at $\zeta\approx 8.084$. The bifurcating branch ``PO3'' 
differs from PO2 by the shorter period at bifurcation, and 
a rather different type of oscillation, indicated by 
the two snapshots at $t=0$ and $t=T/2$. While on PO2 the spots 
oscillate synchronously (the vesicle breathes as a whole), 
on PO3 the spots oscillate individually, and change their distances, 
and $\lam$ and $\spr{\phi}$ have smaller amplitude 
than on PO2. In any case, PO3 again bifurcates 
subcritically and hence is initially unstable, and DNS 
falls back on the octahedral branch. Then PO3 shows a fold 
at $\zeta\approx 7.839$ but does not stabilize there. Instead, 
cDNS from PO3 after the fold still tends to the octahedral 
branch. 
The decahedral orange branch loses stability 
at $\zeta\approx 8.340$ and at $\zeta\approx 5.704$, and then at smaller 
$\zeta$  (see sample B), 
seems to bifurcate from the axisymmetric branch analogous 
to the light orange branch in Fig.\ref{3f5a} (see also the 
$t=0$ plot in Fig.\ref{3f6}(b) for a perturbed sample). However,  
we cannot quite make this connection in the continuation, and 
instead the orange branch develops some uncontrolled branch jumping 
near $\zeta=4.45$. 

To close, we remark that 
without the mesh refinement it was not possible 
to follow the octahedral and decahedral branches to ``large'' $\zeta$, 
i.e., $\zeta>7$, say, due to loss of symmetry and branch jumping during 
the continuation. Thus, the mesh refinement from Fig.\ref{3f6} encourages 
the keeping of the respective symmetries, although in a somewhat artificial 
way.

\renewcommand{\refname}{References}
\taskip
\small
\bibliographystyle{alpha}\bibliography{./hubib.bib}

@article{G13, 
author = {N. Gov}, 
title = {Guided by curvature: shaping cells by coupling curved membrane proteins and cytoskeletal forces}, 
journal = {Phil.Trans.Royal Society B}, 
volume = {373},
pages = {20170115},
year  = {2017}}

@article {SGS07, 
    AUTHOR = {A. Gomez-Marin and J. Garcia-Ojalvo and J. M. Sancho}, 
 TITLE = {Self-Sustained Spatiotemporal Oscillations Induced by Membrane-Bulk Coupling}, 
   JOURNAL = {Phys.~Rev.~Let.}, 
    VOLUME = {98},
      YEAR = {2007},
     PAGES = 168303}

@article {SG11, 
    AUTHOR = {C. Sample and A. Golovin},  
 TITLE = {Morphological and chemical oscillations in a double-membrane
              system},
   JOURNAL = {SIAM J. Appl. Math.},
    VOLUME = {71},
      YEAR = {2011},
    NUMBER = {2},
     PAGES = {622--634}}

@article {GW16, 
    AUTHOR = {J. Gou and M. Ward}, 
 TITLE = {Oscillatory dynamics for a coupled membrane-bulk diffusion
              model with {F}itzhugh-{N}agumo membrane kinetics},
   JOURNAL = {SIAM J. Appl. Math.},
  FJOURNAL = {SIAM Journal on Applied Mathematics},
    VOLUME = {76},
      YEAR = {2016},
    NUMBER = {2},
     PAGES = {776--804}}

@article {ZHL17, 
    AUTHOR = {Zhao, S. and Healey, T. and Li, Q.},
     TITLE = {Direct computation of two-phase icosahedral equilibria of
              lipid bilayer vesicles},
   JOURNAL = {Comput. Methods Appl. Mech. Engrg.},
  FJOURNAL = {Computer Methods in Applied Mechanics and Engineering},
    VOLUME = {314},
      YEAR = {2017},
     PAGES = {164--179}}

@article {DH15, 
    AUTHOR = {Dharmavaram, S. and Healey, T. J.},
     TITLE = {On the equivalence of local and global area-constraint
              formulations for lipid bilayer vesicles},
   JOURNAL = {Z. Angew. Math. Phys.},
  FJOURNAL = {Zeitschrift f\"{u}r Angewandte Mathematik und Physik. ZAMP.
              Journal of Applied Mathematics and Physics. Journal de
              Math\'{e}matiques et de Physique Appliqu\'{e}es},
    VOLUME = {66},
      YEAR = {2015},
    NUMBER = {5},
     PAGES = {2843--2854}}

@article{LDS20, 
  title={Self-Oscillating Membranes: Chemomechanical Sheets Show Autonomous
Periodic Shape Transformation}, 
  author={I. Levin and R. Deegan and E. Sharon}, 
  journal={Phys. Rev. Lett.},
  volume={125},
  pages={178001},
  year={2020}}

@article{LMAO21, 
  title={Chemically controlled pattern formation in
self-oscillating elastic shells}, 
  author={S. Lia and D. Matoz-Fernandeza and A. Aggarwal and M. 
{Olvera de la Cruz}}, 
  journal={PNAS}, 
  volume={118},
  pages={e2025717118},
  year={2021}}

@article{MSD18, 
  title={Geometry of Wave Propagation on Active Deformable Surfaces}, 
  author={P. Miller and N. Stoop and J. Dunkel}, 
  journal={Phys. Rev. Lett.},
  volume={120},
  pages={268001},
  year={2018}}

@article {TUSY14, 
    AUTHOR = {Ryota Tamate and Takeshi Ueki and Mitsuhiro Shibayama 
and Ryo Yoshida}, 
     TITLE = "{Self-Oscillating Vesicles: Spontaneous Cyclic Structural Changes of
Synthetic Diblock Copolymers}", 
JOURNAL = {Angew.Chem}, 
    VOLUME = {53}, 
pages={11248-11252}, 
      YEAR = {2014}}

@article {SZHZ20, 
    AUTHOR = {Qing Shao and Shaodong Zhang and Zhen Hu and Yongfeng Zhou}, 
     TITLE = "{Multimode Self-Oscillating Vesicle Transformers}", 
JOURNAL = {Angew.Chem}, 
    VOLUME = {59}, 
pages={17125 – 17129}, 
      YEAR = {2020}}

@article {Nog25, 
    AUTHOR = {Hiroshi Noguchi}, 
     TITLE = "{Nonequilibrium Membrane Dynamics Induced by Active Protein Interactions and Chemical Reactions: A Review}", 
JOURNAL = {ChemSystemsChem}, 
    VOLUME = {7}, 
pages={e202400042}, 
      YEAR = {2025}}

@article {DP05, 
    AUTHOR = {Dias, A.P.S. and Paiva, R.C.},
     TITLE = {A note on {H}opf bifurcation with dihedral group symmetry},
   JOURNAL = {Glasg. Math. J.},
  FJOURNAL = {Glasgow Mathematical Journal},
    VOLUME = {48},
      YEAR = {2006},
    NUMBER = {1},
     PAGES = {41--51}}

@misc{NDV25, 
author="D. Nesenberend and  A. Doelman and  Fr. Veerman", 
title="Curvature induced patterns: A geometric, analytical approach to understanding a mechanochemical model", 
note="{Preprint}", 
year="2025"}

@article{L22, 
    AUTHOR = {R. Lipowsky}, 
     TITLE = {Remodeling of Membrane Shape and Topology
by Curvature Elasticity and Membrane Tension}, 
     PAGES = {2101020}, 
JOURNAL ={Adv.Biology}, 
volume=6, 
      YEAR = {2022}}

@article {LRM18, 
    AUTHOR = {Th. Litschel and B. Ramm and R. Maas and M. Heymann and
P. Schwille}, 
     TITLE = {Beating Vesicles: Encapsulated Protein Oscillations
Cause Dynamic Membrane Deformations}, 
   JOURNAL = {Angew. Chem. Int. Ed}, 
    VOLUME = {57},
      YEAR = {2018},
     PAGES = {16286--16290}}

@article {CLSL21, 
    AUTHOR = {S. Christ and Th. Litschel and P. Schwille and R. Lipowsky}, 
     TITLE = {Active shape oscillations of giant vesicles with cyclic closure and opening of membrane necks}, 
   JOURNAL = {Soft Matter}, 
    VOLUME = {17},
      YEAR = {2021},
     PAGES = {319-330}}

@article{YFB22, 
    AUTHOR = {A. Yochelis and S. Flemming and C. Beta}, 
     TITLE = {Versatile patterns in the actin cortex of motile cells:
{S}elf-organized pulses can coexist with macropinocytic ring-shaped waves}, 
JOURNAL = {Phys.Rev.Letters}, 
    VOLUME = {129},
      YEAR = {2022},
     PAGES = {088101}}

@article {BMN05, 
    AUTHOR = {B\"ansch, E. and Morin, P. and Nochetto, R.},
     TITLE = {A finite element method for surface diffusion: the parametric
              case},
   JOURNAL = {J. Comput. Phys.},
  FJOURNAL = {J. of Computational Physics},
    VOLUME = {203},
      YEAR = {2005},
    NUMBER = {1},
     PAGES = {321--343}}

@article{GNPS96,  
author = {R.E. Goldstein and P. Nelson and T. Powers and U. Seifert}, 
title = {Front progagation in the pearling instability of tubular vesicles}, 
journal = {J.Phys.II}, 
volume = {6},
pages = {767--796},
year  = {1996}}

@article {NY12, 
    AUTHOR = {Nagasawa, Takeyuki and Yi, Taekyung},
     TITLE = {Local existence and uniqueness for the {$n$}-dimensional
              {H}elfrich flow as a projected gradient flow},
   JOURNAL = {Hokkaido Math. J.},
  FJOURNAL = {Hokkaido Mathematical Journal},
    VOLUME = {41},
      YEAR = {2012},
    NUMBER = {2},
     PAGES = {209--226},
       DOI = {10.14492/hokmj/1340714413},
       URL = {https://doi.org/10.14492/hokmj/1340714413}}

@misc{rupp24, 
    AUTHOR = {F. Rupp and C. Scharrer and M. Schlierf}, 
     TITLE = "{Gradient flow dynamics for cell membranes in the Canham-Helfrich model, \url{https://arxiv.org/pdf/2408.07493}}", 
      YEAR = {2024}}

@misc{SGJW22, 
    AUTHOR = {O. Stein and E. Grinspun and A. Jacobson and Wardetzky, M.},
     TITLE = {A mixed finite element method with piecewise linear elements for the biharmonic equation on surfaces, \url{http://arxiv.org/pdf/1911.08029}},
      YEAR = {2022}}

@article{IMKYG13, 
author = {E. B. Isaac and  U. Manor and B. Kachar and A. Yochelis and  N. Gov}, 
title = {Linking actin networks and cell membrane via a reaction-diffusion-elastic description of nonlinear filopodia initiation}, 
journal = {Phys.Rev.E},
volume = {88},
pages = {022718-1},
year  = {2013}}

@article {MHMC13,
    doi = {10.1371/journal.pone.0082617},
    author = {Mercker, M. AND Hartmann, D. AND Marciniak-Czochra, A.},
    journal = {PLOS ONE},
    publisher = {Public Library of Science},
    title = {A Mechanochemical Model for Embryonic Pattern Formation: Coupling Tissue Mechanics and Morphogen Expression},
    year = {2013},
    month = {12},
    volume = {8},
    url = {https://doi.org/10.1371/journal.pone.0082617},
    pages = {1-6}}

@article {MMRH13, 
    AUTHOR = {Mercker, M. and Marciniak-Czochra, A. and Richter, T. and
              Hartmann, D.},
     TITLE = {Modeling and computing of deformation dynamics of
              inhomogeneous biological surfaces},
JOURNAL = {SIAM J. Appl. Math.},
    VOLUME = {73},
      YEAR = {2013},
    NUMBER = {5},
     PAGES = {1768--1792}}

@article {BMRM18, 
    AUTHOR = {Brinkmann, F. and Mercker, M. and Richter, T. and Marciniak-Czochra, A.}, 
     TITLE = {Post-{T}uring tissue pattern formation: {A}dvent of mechanochemistry},
   JOURNAL = {PLOS Computational Biology}, 
      YEAR = {2018},
     PAGES = {1006259}}

@article {TN20, 
    AUTHOR = {Naoki Tamemoto and Hiroshi Noguchi}, 
     TITLE = {Pattern formation in reaction–diffusion system on membrane with mechanochemical feedback}, 
   JOURNAL = {Scientific reports}, 
    VOLUME = {10},
      YEAR = {2020},
    NUMBER = {19582}}

@article {TN21, 
    AUTHOR = {Naoki Tamemoto and Hiroshi Noguchi}, 
     TITLE = {Reaction-diffusion waves coupled with membrane curvature}, 
   JOURNAL = {Soft Matter}, 
    VOLUME = {17},
      YEAR = {2021},
    NUMBER = {6589--6596}}

@article {TN22, 
    AUTHOR = {Naoki Tamemoto and Hiroshi Noguchi}, 
     TITLE = {Excitable reaction-diffusion waves of curvature-inducing proteins on deformable membrane tubes}, 
   JOURNAL = {Phys.Rev.E}, 
    VOLUME = {106},
      YEAR = {2022},
    NUMBER = {024403}}

@misc{mcSI,
    AUTHOR = {Meiners, A. and Uecker, H.},
     TITLE = "{Supplementary information for {\em 
Breathing and moving vesicles in a geometric mechanochemical model}, 
\url{https://pde2path.uol.de/apps}}", 
      YEAR = {2026}}

@misc{geomtut,
    AUTHOR = {Meiners, A. and Uecker, H.},
     TITLE = {Differential 
 geometric bifurcation problems in {\tt pde2path} -- algorithms and tutorial 
examples}, 
note="Available at \cite{p2phome}", 
      YEAR = {2024}}

@article{mucyl, 
    AUTHOR = {Meiners, A. and Uecker, H.},
     TITLE = {Helfrich cylinders -- instabilities, bifurcation analysis and 
amplitude equations}, 
journal="SIAM J. Appl. Dyn.Systems", 
      YEAR = {2026}}

@misc{AMPHD, 
    AUTHOR = {Meiners, A.}, 
     TITLE = {Differential geometric bifurcation problems -- Theory and applications to minimal surfaces and biomembranes}, 
note="PhD thesis", 
      YEAR = {2024}}

@article{geompap,
    AUTHOR = {Meiners, A. and Uecker, H.},
     TITLE = {Numerical continuation and bifurcation for 
 differential geometric {PDEs}}, 
journal="{Numerical Mathematics -- Theory Methods Applications}", 
volume={OA-2024}, 
pages={0005}, 
      YEAR = {2024}}

@article{KPP17,
	author = {M. Koiso and B. Palmer and P. Piccione},
	journal = {Journal of the Mathematical Society of Japan},
	number = {4},
	pages = {1519 -- 1554},
	title = {Stability and bifurcation for surfaces with constant mean curvature},
	url = {https://doi.org/10.2969/jmsj/06941519},
	volume = {69},
	year = {2017},
	}

@article {H73, 
    AUTHOR = {Helfrich, W.},
     TITLE = {Elastic properties of lipid bilayers: Theory and possible experiments},
   JOURNAL = {Zeitschrift f\"ur Naturforschung}, 
    VOLUME = {28},
      YEAR = {1973},
     PAGES = {693}}

@misc{gpgit, 
    AUTHOR = {A. Jacobson}, 
     TITLE = {gptoolbox, \url{https://github.com/alecjacobson/gptoolbox}}, 
      YEAR = {2024}}

@article{sei97, 
author = {U.  Seifert },
title = {Configurations of fluid membranes and vesicles},
journal = {Advances in Physics},
volume = {46},
number = {1},
pages = {13-137},
year  = {1997},
publisher = {Taylor & Francis}}

@incollection{SL95, 
author = {U.  Seifert and R. Lipowsky},
title = {Morphology of {V}esicles}, 
 BookTitle = {Handbook of Biological Physics},
volume=1, 
editor={R. Lipowsky and E. Sackmann}, 
pages = {403--463},
year  = {1995},
publisher = {Elsevier}}

@article{SBL90, 
author = {U.  Seifert and K. Berndl and R. Lipowsky},
title = {Shape transformations of vesicles: {P}hase diagram },
journal = {Phys.Rev. A}, 
volume = {44},
pages = {1182--1202},
year  = {1990}}

@article {PXP20, 
    AUTHOR = {Paquin-Lefebvre, Fr. and Xu, Bin and DiPietro, K.L.
              and Lindsay, A.E. and Jilkine, A.},
     TITLE = {Pattern formation in a coupled membrane-bulk
              reaction-diffusion model for intracellular polarization and
              oscillations},
   JOURNAL = {J. Theoret. Biol.},
  FJOURNAL = {Journal of Theoretical Biology},
    VOLUME = {497},
      YEAR = {2020},
     PAGES = {110242, 23}}

@article{BEGY23, 
    AUTHOR = {C. Beta and L. {Edelstein--Keshet} and N. Gov and  A. Yochelis}, 
     TITLE = {From actin waves to mechanism and back: How theory aids 
biological understanding}, 
JOURNAL = {Elife}, 
    VOLUME = {12},
      YEAR = {2023},
     PAGES = {e87181}}

@article{HMY25, 
    AUTHOR = {J.M. Hughes and C. Martinez--Torres and C. Beta and L. {Edelstein--Keshet} and  A. Yochelis}, 
     TITLE = {A dissipative mass conserved reaction–diffusion
system reveals switching between coexisting
polar and oscillatory cell motility states}, 
JOURNAL = {Chaos}, 
    VOLUME = {35},
      YEAR = {2025},
     PAGES = {051103}}

@book {p2pbook, 
    AUTHOR = {Uecker, H.},
     TITLE = {Numerical continuation and bifurcation in Nonlinear PDEs}, 
 PUBLISHER = {SIAM, Philadelphia, PA},
      YEAR = {2021}}

@misc{p2phome, 
author="pde2path",  
title="{\url{https://pde2path.uol.de/}}", 
year="2026"}

@incollection {Matt04, 
    AUTHOR = {Matthews, P. C.},
     TITLE = {Pattern formation on a sphere},
 BOOKTITLE = {Dynamics and bifurcation of patterns in dissipative systems},
    SERIES = {World Sci. Ser. Nonlinear Sci. Ser. B},
    VOLUME = {12},
     PAGES = {102--123},
 PUBLISHER = {World Sci. Publ., Hackensack, NJ},
      YEAR = {2004}}

@book{hlr89, 
author="Hairer, E. and Lubich, Ch. and Roche, M.",
title="{The Numerical Solution of Differential-Algebraic Systems by Runge-Kutta Methods}",
publisher="{Springer}",
year="1989"}

@book{hw2,
author="Hairer, E. and Wanner, G.",
title="{Solving ordinary differential equations. II: Stiff and
    differential-algebraic problems.}",
publisher="{Springer}",
year="1996"}

@book{ae1,
author="Amann, H. and Escher, J.",
title="{Analysis I.}",
publisher="{Birkh\"auser}",
year="2006"}

@book{ae2,
author="Amann, H. and  Escher, J.",
title="{Analysis II.}",
publisher="{Birkh\"auser}",
year="2006"}

@book{ae3,
author="Amann, H. and Escher, J.",
title="{Analysis III.}",
publisher="{Birkh\"auser}",
year="2008"}
\end{document}